\documentclass[10pt]{amsart}
\usepackage{amssymb}
\usepackage{MnSymbol}

\usepackage{amscd}
\usepackage{amsthm}
\usepackage[dvips]{graphicx}
\usepackage{color}
\usepackage[all]{xy}
\usepackage{verbatim}
\usepackage{cite}

\newtheorem{Lemma1}{{Lemma}}[section]
\newtheorem{Theo1}[Lemma1]{{Theorem}}
\newtheorem*{Theo2}{{Theorem}}
\newtheorem{Def1}[Lemma1]{{Definition}}
\newtheorem{Prop1}[Lemma1]{{Proposition}}
\newtheorem{Claim1}[Lemma1]{{Claim}}
\newtheorem{Rem1}[Lemma1]{{Remark}}
\newtheorem{Cor1}[Lemma1]{{Corollary}}
\newtheorem{Ex1}[Lemma1]{{Example}}
\newtheorem{Not1}[Lemma1]{{Notation}}

\newenvironment{Lemma}{\begin{Lemma1}}{\end{Lemma1}}
\newenvironment{Def}{\begin{Def1}\rm}{\end{Def1}}
\newenvironment{Prop}{\begin{Prop1}}{\end{Prop1}}
\newenvironment{Rem}{\begin{Rem1}\rm}{\end{Rem1}}
\newenvironment{Theorem}{\begin{Theo1}}{\end{Theo1}}

\newenvironment{Remark}{\begin{Rem1}\rm}{\end{Rem1}}
\newenvironment{Cor}{\begin{Cor1}}{\end{Cor1}}

\title{Differential graded Brauer groups over dg-rings}

\author{Xiaoxiao Xu}
\address{\newline
X.X.: Universit\'e de Picardie,
\newline D\'epartement de Math\'ematiques et LAMFA (UMR 7352 du CNRS),
\newline 33 rue St Leu,
\newline F-80039 Amiens Cedex 1,
\newline France\newline 
and\newline
School of Mathematical Sciences,\newline  
Shanghai Key Laboratory of PMMP,\newline 
East China
Normal University, \newline
Shanghai 200241,\newline China}
\email{xiaoxiao.xu@u-picardie.fr}

\author{Alexander Zimmermann}
\address{\newline
A.Z.: Universit\'e de Picardie,
\newline D\'epartement de Math\'ematiques et LAMFA (UMR 7352 du CNRS),
\newline 33 rue St Leu,
\newline F-80039 Amiens Cedex 1,
\newline France}
\email{alexander.zimmermann@u-picardie.fr}

\date{June 20, 2026}

\newcommand{\lra}{\longrightarrow}

\newcommand{\ra}{\rightarrow}
\newcommand{\sdp}{\rtimes}%{\times\kern-.2em\vrule height1.1ex depth-.05ex}
\newcommand{\epi}{\lra \kern-.8em\ra}

\newcommand{\Z}{{\mathbb Z}}

\newcommand{\im}{\textup{im}}
\newcommand{\dickebox}{{\vrule height5pt width5pt depth0pt}}

\newcommand{\id}{\textup{id}}

\newcommand{\Hom}{\textup{Hom}}

\newcommand{\End}{\textup{End}}

\newcommand{\Br}{\textup{Br}}
\newcommand{\dgBr}{\textup{dgBr}}
\newcommand{\grBr}{{\textup{grBrW}}}
\newcommand{\ogrBr}{\textup{grBr}}

 \normalbaselines

\begin{document}

\begin{abstract}
Brauer groups of graded rings were defined and studied by Caenepeel and van Oystaeyen.
We study the question what happens to this group if a ${\mathbb Z}$-graded ring carries 
a differential graded structure in addition. 
We define a Brauer group for differential graded algebras over differential graded 
graded-commutative or commutative base rings. Based on previous work we give an 
explicit classification of dg-fields, and compute as an example 
the so-defined Brauer group in each case
explicitly. 
\end{abstract}

\maketitle

\section*{Introduction}

Differential graded algebras (dg-algebras for short)
were defined by Cartan in 1954 \cite{Cartandg}, 
and proved to be a most essential tool in homological algebra, 
algebraic topology, algebraic geometry,
differential geometry and many other subjects. Most astonishing, 
a ring theoretic point of view 
was not considered, until very recently. Aldrich and Garcia-Rozas 
characterized acyclic dg-algebras completely \cite{Tempest-Garcia-Rochas}.
Orlov \cite{Orlov1} 
considered finite dimensional dg-algebras over a field from an algebraic 
geometry point of view. Goodbody \cite{Goodbody}, 
the second author \cite{dgorders,dgBrauer,dgGoldie,dgfields,dgseparable,dgHopkins} 
and Orlov \cite{Orlov2} studied 
the ring theory of dg-algebras with respect to several point of views and directions. 
The purpose of this present paper is to continue our aim to understand
the ring theoretic behaviour of dg-algebras.
 
A starting point for the present paper is that 
in \cite{dgfields} dg-division algebras were defined as those
dg-algebras without non trivial left or right dg-ideal, and that it could be 
shown that this is 
equivalent to the property that the cycles of the dg-algebras 
form a graded-division algebra
(cf \cite{vanGeel}, see also
\cite{gradedrings,NastasescuVanOystaen}). It is proved there that
dg-division algebras are either acyclic or the differential is $0$ and then 
the algebra is a graded-division ring.  
This is another case of our general observation is that 
frequently the property needed 
for the dg-algebra is 
asked for the graded ring formed by the cycles of the dg-algebra.  

\subsection*{On Brauer groups, a historical summary.}
The Brauer group $\textup{Br}(K)$ of a field $K$ parameterizes central simple $K$-algebras, and it is 
a fundamental classical result that $\textup{Br}(K)$ is isomorphic to the degree $2$ group cohomology
of the Galois group of the separable closure $K^{sep}$ over $K$ with values in the multiplicative group
of $K^{sep}$. C.T.C.~Wall \cite{CTCWall} then considered algebras graded by a cyclic 
group of order $2$, defined and studied the then commonly called Brauer-Wall group analogously 
with a tensor product obeying the Koszul sign rule.  
Auslander and Goldman~\cite{AuslanderGoldman} generalised the classical ungraded 
situation to the so-called
Brauer-Azumaya group of a commutative ring $R$, and Small~\cite{CharlesSmall} defined and studied  
the Brauer-Wall group for a commutative ring. Its elements are Morita equivalence classes of Azumaya 
algebras over this base ring, where an $R$ algebra is Azumaya if it is finitely generated 
projective over $R$ and $A\otimes_RA^{op}$ is isomorphic to $\End_R(A)$ in a very precise way (see 
Section~\ref{versionsofbrauergroups} for a precise definition). 
Caenepeel \cite{CaenepeelBrauergroupsHopfalgebras}, as well as Caenepeel and van Oystaeyen 
\cite{caenepeelvanoystaeyen}
then considered rings graded by a  group and defined and studied 
graded Azumaya algebras. 
Grothendieck~\cite{GrothgroupedeBrauerI,GrothgroupedeBrauerII,GrothgroupedeBrauerIII} 
considered schemes of algebras 
and proved that under favorable hypotheses the Brauer group of the scheme $X$ of algebras 
embeds into the torsion part of the degree $2$ \'etale cohomology of the scheme with 
values in ${\mathbb G}_m$. For a detailed account on the Brauer-Grothendieck group we refer to 
\cite{Colliotthelenes}.

\subsection*{What happens if a differential graded structure is present. Outline of our results.}
In the present paper we answer the question what happens to the graded Brauer group of 
a $\Z$-graded algebra in the sense of \cite{CaenepeelBrauergroupsHopfalgebras,caenepeelvanoystaeyen}
when the graded Azumaya algebras carry in addition a structure
of a differential graded algebra. 
More generally, we answer the question how to define and study a Brauer group of a 
commutative or graded-commutative dg-algebra. We need to answer first the question what 
should be a dg-Azumaya algebra. 

We say that a dg-algebra $(A,d_A)$ over a graded-commutative 
dg-algebra $(K,d_K)$ is dg-Azumaya of the second kind 
if $A$ is graded-Azumaya over $K$. We prove in 
Proposition~\ref{groupstructureofdgbrauerofsecondkind}
that this gives a well-defined dgBrauer group
denoted $\dgBr^{II}(K,d_K)$. 
We should observe that a graded-commutative algebra
$A=\bigoplus_{n\in\Z}A_n$  graded over $\Z$, is actually 
a graded algebra with group being cyclic of order $2$. 
$$A=\bigoplus_{n\in\Z}A_n=\bigoplus_{n\in\Z}A_{2n}\oplus \bigoplus_{n\in\Z}A_{2n+1}
=B^{(0)}\oplus B^{(1)}$$
where $B^{(0)}=\bigoplus_{n\in\Z}A_{2n}$ and $B^{(1)}=\bigoplus_{n\in\Z}A_{2n+1}$.
Though similar, this situation is not quite the situation considered by C.T.C. Wall, since 
each of the components $B^{(0)}$ and $B^{(1)}$ is graded over $\Z$, whereas C.T.C.~Wall considers
the case when $B^{(0)}$ and $B^{(1)}$ are ungraded. 
We hence get that $\dgBr^{II}(K,d_K)$ combines aspects of the Brauer-Wall group, and 
aspects of Caenepeel's graded-Brauer group. We indicate in Section~\ref{gradedcommutativegradedSection}
their relation. 

The experience in our previous investigations on 
the ring theory of dg-algebras indicates that the graded ring of cycles should control the 
behaviour. This motivates the alternative and more direct definition of a dg-Brauer group. 
A dg-algebra $(A,d_A)$ over a (graded)-commutative dg-base ring $(K,d_K)$ is dg-Azumaya 
of the first kind if $\ker(d_A)$ is graded-Azumaya over $\ker(d_K)$ 
(cf \cite{AuslanderGoldman, caenepeelvanoystaeyen}). Equivalence classes are
those given by the cycles being equivalent in the graded setting.
The product of two such algebras is the equivalence class of the tensor product of the 
cycles over $\ker(d_K)$ and then tensored over $\ker(d_K)$ with $K$, giving by this the 
differential obtained from $K$. It turns out that his defines 
a well-defined dgBrauer group 
(cf Proposition~\ref{grouplawindgbrI}), called 
$\dgBr^I(K,d_K)$. It is not hard to show that this group 
is isomorphic to the graded Brauer-Wall group of $\ker(d_K)$ 
(cf Lemma~\ref{dgBrIisgrBr} and Appendix~\ref{Appendixsection}). 

We always have a homomorphism 
$$\dgBr^{II}(K,d_K)\lra \dgBr^I(K,d_K)=\ogrBr(\ker(d_K))$$
which proves to be an isomorphism in case of  
$(K,d_K)$ having semisimple dg-module category 
(cf Theorem~\ref{theultimateforacyclic}), and proves 
to be a split epimorphism in case of $d_K=0$ (cf Theorem~\ref{dgBrIIvsusgradedBrK}). 
If $(K,d_K)$ is more general, the group $\dgBr^{II}(K,d_K)$ is related to bad gradings 
on matrix algebras~\cite{DascalescuIonNastasescuRioMontes}, and seems to encode 
sophisticated invariants of a ring. 
 
We are able to compute the dgBrauer groups of both types in case of dg-fields $(K,d_K)$ 
in each case completely (cf Section~\ref{Examplesection}).

\subsection*{Quadratic forms}
Suppose that multiplication by $2$ is an injective map on $K$. Then
the Brauer-Wall group is closely linked to the classification of quadratic 
forms via a homomorphism of the Grothendieck-Witt group to the Brauer-Wall group,
associating to a non degenerate bilinear form $b$ on a faithful projective $K$-module $V$ 
its Cifford algebra, presented as a quotient of the tensor algebra on $V$ as 
$$\textup{Cliff}_b(V):=TV/(x\otimes y+y\otimes x-2\cdot b(x,y);\; x,y\in V\textup{ homogeneous}).$$ 
This algebra then gives a representative 
of a $\Z/{2\Z}$-graded-Azumaya algebra in the Brauer-Wall group. One might ask 
if this is possible for our dg-version. Generalizing slightly Chen and Kang  
\cite[Examples 3.3]{ChenKang}, we get for a field $K$ concentrated in degree $0$, and $V$ 
a $2n$-dimensional $\Z$-graded vector space with basis $x_1,\dots,x_n$ in degree $1$ and $y_1,\dots,y_n$ in degree $-1$ with skew symmetric bilinear form defined by 
$b(x_i,y_j):=\delta_{i,j}$, the Kronecker symbol, then the 
$\Z$-graded Clifford algebra 
$$\textup{Cliff}_b(V):=TV/(x\otimes y+(-1)^{|x||y|}y\otimes x-2\cdot b(x,y);\; x,y\in V\textup{ homogeneous})$$ is isomorphic 
to the Weyl algebra $W_n(K)$, with graded centre $K$, 
but is infinite dimensional.  
Hence the $\Z$-graded Clifford algebra of a graded module 
is not a graded-Azumaya algebra, which shows that the link between 
quadratic spaces and Brauer-Wall groups does not generalize to the dg-setting.

\subsection*{Ruminations on the derived version.}
To\"en defined in \cite{ToenAzumaya} an alternative concept, not using 
the classical Morita type equivalence relation, but asking instead that the isomorphism 
of $A\otimes_KA^{op}$ with $\End_K(A)$ is a quasi-isomorphism and the equivalence relation alike. 
This then produces a very different concept. In particular, since dg-fields are mostly acyclic, 
To\"en's construction vanishes on these dg-fields. 
Vezzosi \cite{Vezzosi} produced a Clifford algebra and a Grothendieck-Witt 
paralleling To\"en's version of derived Azumaya algebras. However, as is mentioned in 
\cite[Remark 2.]{Vezzosi} the correspondence to To\"en's work does not quite fit. 
Further, again Vezzosi's construction is invariant under derived equivalences, whence vanishes in case of acyclic dg-base fields. 

By the above mentioned reason, this approach is unsuitable for our purpose.

\subsection*{Our paper is organized as follows. }
In Section~\ref{notationsection} we first recall some basic facts about dg-algebras and
dg-modules, also in order to fix our conventions and notations. Then, we add some 
concepts for the constructions over (graded) commutative dg-algebras. Further, we recall
our results on dg-division algebras. As a complement to our results in \cite{dgfields}
as our first main result Theorem~\ref{fieldclassscases}, 
we classify completely  in Section~\ref{dgfieldsclassification} 
commutative and graded-commutative dg-division algebras.  
In Section~\ref{definingdgbrauersection} we define our two different 
concepts of dgBrauer groups, and show that these provide groups in each case. 
In Section~\ref{differentialzerosection} we show in
Theorem~\ref{dgBrIIvsusgradedBrK} 
that the dgBrauer group of the second kind maps always to the 
dgBrauer group of the first kind and this map is split epic if 
the differential of the base ring is $0$. This is our second main result. 
Further, for acyclic base ring we prove that there is an isomorphism of the 
two concepts, 
which presents our second main result Theorem~\ref{theultimateforacyclic}.
Finally, we provide Examples in Section~\ref{Examplesect}. In particular, 
we compute the dg-Brauer groups of each type explicitly and completely 
for each of the cases of base dg-fields from Theorem~\ref{fieldclassscases}.
In Appendix~\ref{Appendixsection} 
we show that for a graded algebra over a graded commutative ring the 
equivalence classes of graded Azumaya algebras form a graded Brauer group. 
This may be well-known to the experts, but we 
could not find an explicit treatment in the literature.

\subsection*{Acknowledgement:}
After the second author lectured on \cite{dgorders,dgBrauer} during 
the conference 'journ\'ees d'alg\`ebres' 
in June 2023 at the universit\'e Blaise Pascal in Clermont-Ferrand, 
Bernhard Keller suggested that 
\cite{dgBrauer} should be considered for (differential) graded base rings. 
By email from May 5, 2025 
Ziqi Liu asked the second author if a result like the one in our manuscript 
could be possible, in the perspective of studying the dg Brauer 
groups of the dg-category 
of a (good) scheme. After the second author lectured on \cite{dgfields} 
on the occasion of 
the 'Mini-Workshop on Algebra' which took 
place in September 2025 at the University of Science and Technology of China 
in Hefei, Roozbeh Hazrat insisted 
that a result like the one in the present paper should be explored.  
We are very grateful to Bernhard Keller, and to Roozbeh Hazrat for 
their suggestion. Further, we thank Ziqi Liu for helpful comments 
improving an earlier version of the paper, and
we thank Pierre Baumann for sending us a copy of \cite{TilborghsVanOystaeyen}.
Finally, we are grateful to the Chinese Scholarship Council for their generous support 
with a scholarship  
making our collaboration possible.

\section{Notations, conventions, and generalities on differential graded algebras}
\label{notationsection}

\subsection{Elementary definitions of dg-algebras and modules} 

Let $K$ be a commutative base ring. A 
differential graded algebra (dg-algebra for short) over $K$ is a 
$\Z$-graded $K$-algebra $A$ which is equipped with a $K$-linear 
map $d:A\lra A$ which is homogeneous of degree $1$, with $d^2=0$ and 
such that for all homogeneous $a,b\in A$ we get 
$$d(a\cdot b)=d(a)\cdot b+(-1)^{|a|}a\cdot d( b),$$
where we denote by $|b|$ the degree of $b\in A$. 
The opposite algebra $A^{op}$ is $A$ with the same additive structure as
$A$ and multiplication denote by $\cdot_{op}$ and defined as
$$b\cdot_{op}a:=(-1)^{|a||b|}a\cdot b$$ 
for all homogeneous $a,b\in A$. Then $(A^{op},d)$ is again a dg-algebra, 
with the same differential $d$.

A left dg-module $(M,\delta)$ over a dg-algebra $(A,d)$ over $K$ is a 
$\Z$-graded left $A$-module with a $K$-linear endomorphism $\delta:M\lra M$ of 
degree $1$ with $\delta^2=0$ and  
$$\delta(a\cdot m)=d(a)\cdot m+(-1)^{|a|}a\cdot \delta(m)$$
for all homogeneous $a\in A$ and $m\in M$. 
A right dg-module $(M,\delta)$ over $(A,d)$ is a left dg-module over $(A^{op},d)$.
A homomorphism of dg-algebras $\varphi:(A,d_A)\lra (B,d_B)$
is a graded degree $0$ homomorphism of algebras $\varphi:A\lra B$ such that 
$\varphi\circ d_A=d_B\circ\varphi$. 

If $(M,\delta_M)$ is a left dg-module over $(A,d_A)$, and $(N,\delta_N)$ is a 
right dg-module over $(A,d_A)$, then 
$N\otimes_AM$ allows a differential $\delta_{N\otimes_AM}$
given by 
$$\delta_{N\otimes_AM}(n\otimes m):=\delta_N(n)\otimes m+(-1)^{|n|}n\otimes\delta_M(m)$$
for all homogeneous $n\in N$ and $m\in M$.

If $(M,\delta_M)$ and $(N,\delta_N)$ are left dg-modules over $(A,d_A)$. 
Then 
$$\Hom_A^\bullet(M,N):=
\bigoplus_{k\in\Z}\Hom_A^k(M,N)$$
and 
$$\Hom_A^k(M,N):=\left\{f:M\lra N\;|\;\begin{minipage}{8.5cm}\begin{center} $
\forall_{\ell\in\Z} f(M_\ell)\subseteq N_{\ell+k}\textup{ and }$\\
$\forall_{a\in A,m\in M\textup{ homogeneous,}}\;f(am)=(-1)^{|a|k}af(m)\;\textup{ and}$\\
$f\textup{ additive}$
\end{center}\end{minipage}\right\}$$
$\Hom_A^\bullet(M,N)$ is equipped with a differential 
$$d_{\Hom}(f):=\delta_N\circ f-(-1)^{|f|}f\circ\delta_M.$$
Denote $\End_A^\bullet(M,\delta_M):=\Hom_A^\bullet(M,M)$, and likewise for $N$, then 
$(\End_A^\bullet(M,\delta_M),d_\Hom)$ is a dg-algebra and likewise for $(N,\delta_N)$. 
Further, $(\Hom_A^\bullet(M,N),d_{\Hom})$ is a differential graded bimodule over
$(\End_A^\bullet(N,\delta_N),d_\Hom)-(\End_A^\bullet(M,\delta_M),d_\Hom)$. 

A homomorphism of dg-modules 
$(M,\delta_M)\lra (N,\delta_N)$ over the same dg-algebra $(A,d_A)$
is an element in the degree $0$ cycles of $(\Hom^\bullet_A((M,\delta_M),(N,\delta_N)),d_\Hom)$. 

Let $(A,d_A)$ be a dg-algebra. Then the left regular module is a left dg-module
over $(A,d_A)$, and likewise the 
right regular module is a right dg-module over $(A,d_A)$. Also, $(A,d_A)$ is a dg-bimodule
over $(A,d_A)-(A,d_A)$. A dg left ideal of $(A,d_A)$ is a left dg-submodule of 
the left regular module, and likewise for right and twosided dg-ideals.

For a dg-algebra $(A,d)$ we denote by $dg-mod(A,d)$ the category of 
finitely generated left dg-modules over $A$,
and by $gr-mod(A)$ the category of finitely generated graded left modules over $A$. 

\subsection{dg-algebras over (graded) commutative dg-algebras}
If $(R,d_R)$ is a graded-commutative dg-algebra over $K$, then we shall define 
the structure of a dg-algebra $(A,d_A)$ over 
the dg-algebra $(R,d_R)$. 

Recall that for a commutative ring $K$ a $K$-algebra $A$ is a ring $A$ together with 
a ring homomorphism $$\lambda: K\lra Z(A).$$
If $K$ is graded commutative, this definition does not make sense anymore. However we can 
obtain a meaningful correction easily. 
Let $K$ be a graded-commutative ring. Then a graded $K$-algebra $A$ is a graded 
ring $A$ together with 
a graded ring homomorphism $$\lambda: K\lra Z_{gr}(A).$$

For dg-algebras, the second condition is much more natural for the following reason. 
If $(K,d_K)$ is a dg-algebra itself, then we have
$$d(xy)=d(x)y+(-1)^{|x|}xd(y).$$ 
If $K$ is commutative, then 
$$(-1)^{|y|}d(x)y+xd(y)=d(y)x+(-1)^{|y|}yd(x)=d(yx)=d(xy)=d(x)y+(-1)^{|x|}xd(y).$$ 
Hence 
$$d(x)y(1-(-1)^{|y|})=xd(y)(1-(-1)^{|x|})$$ 
for all homogeneous $x$ and $y$.
If $K$ is graded commutative, then 
$$d(xy)=d(x)y+(-1)^{|x|}xd(y)=
(-1)^{|x||y|}d(y)x+(-1)^{|y|(|x|+1)}yd(x)=(-1)^{|x||y|}d(yx).$$ 
We see that dg-algebras over a commutative dg-ring impose conditions on the base ring,
whereas dg-algebras $A$ over graded-commutative dg-rings do not have to satisfy any 
additional condition.

If $(R,d_R)$ is graded-commutative, and $(A,\delta_A)$ as well as $(B,\delta_B)$ are 
dg-algebras over $(R,d_R)$, then the tensor product $A\otimes_RB$ carries a non standard algebra
structure, and actually 
$(A\otimes_RB,d_{A\otimes_RB})$ is a dg-algebra again by 
$$(a_1\otimes b_1)\cdot(a_2\otimes b_2)=(-1)^{|b_1||a_2|}(a_1a_2\otimes b_1b_2)$$
for all homogeneous $a_1,a_2\in A$ and $b_1,b_2\in B$. 
In order to distinguish this twisted algebra structure we denote it by 
$A\widehat\otimes_RB$. 
This is well-defined since 
$(A,d_A)$ is a dg-algebra over $(R,d_R)$ and hence satisfies the Leibniz rule
$$d_A(a\cdot r)=d_A(a)\cdot r+(-1)^{|a|}a\cdot d_R(r)$$
for all homogeneous $a\in A$ and $r\in R$, and likewise for $(B,d_B)$.
In particular, we cannot forget the differential on $R$, and
in general $(R,d_R)$ is not a dg-algebra over $(R,0)$.  

We should note that the existence of a differential of the tensor product 
of the two algebras implied 
the algebra structure on $A\widehat\otimes_RB$ above, and this in turn implies that
$R$ maps to the graded centre of $A$ and of $B$, and hence
that $R$ is graded-commutative.
Indeed, consider $a,a'\in A$, $b,b'\in B$ and $r\in R$ all homogeneous. 
Then
\begin{eqnarray*}
(ar\otimes b)\cdot(a'\otimes b')\lefteqn{-(a\otimes rb)\cdot(a'\otimes b')}\\
&=&(-1)^{|b||a'|}(ara'\otimes bb')-(-1)^{(|b|+|r|)|a'|}(aa'\otimes rbb')\\
&=&(-1)^{|b||a'|}(ara'\otimes bb')-(-1)^{(|b|+|r|)|a'|}(aa'\otimes rbb')\\
&=&(-1)^{|b||a'|}(a(ra'-(-1)^{|r||a'|}a'r)\otimes bb')
\end{eqnarray*}
which implies $ra'=(-1)^{|r||a'|}a'r$. 
Hence, once we want to form dg-Brauer groups of any form we will need to work with the 
twisted tensor product on algebras as above, and graded-commutative base rings $R$.

Since we could not find a complete enough treatment 
of the theory of graded algebras over 
a graded-commutative algebra 
we reprove in Appendix~\ref{Appendixsection} the necessary ring theoretical background 
in particular the aspects 
necessary to prove that equivalence classes of graded Azumaya algebras 
form a group, the graded Brauer group. 

As we denoted the non standard algebra structure on $A\otimes_RB$ by $A\widehat\otimes_RB$
we recall that if $C$ is a $\Z$-graded algebra, and $M$, $N$ are 
two $\Z$-graded $A$-modules we have two different concepts of graded homomorphism spaces.
We call $\Hom_{C,graded}(M,N)$ the graded space of graded $C$-linear homomorphisms $M\lra N$.
For the dg-hom-complex we used another concept. 
We call $\widehat\Hom_{C,graded}(M,N)$ the graded space generated by graded maps 
$\varphi:M\lra N$ such that $\varphi(cm)=(-1)^{|\varphi||c|}c\varphi(m)$, i.e. 
taking into account the Koszul sign rule.

\subsection{Various versions of Brauer groups}
\label{versionsofbrauergroups}

The most classical Brauer group of a field is defined as follows 
(cf e.g. \cite{MaximalOrders}). 
A finite dimensional $K$-algebra is central simple if it 
is a simple algebra with centre $K$. 
A consequence of the property of being central simple is that 
$A\otimes_KA^{op}$ is Morita equivalent to $K$. The Morita 
equivalence classes of central simple $K$-algebras then form an abelian group, 
the Brauer group $\Br(K)$, with group law given by $-\otimes_K-$.

\subsubsection{The ungraded case:}
Following Azumaya~\cite{Azumaya}, for general commutative rings $R$
one defines an $R$-algebra to be an
Azumaya algebra if 
$A$ is a finitely generated projective generator in the category of $R$-modules, and 
\begin{eqnarray*}
A\otimes_RA^{op}&\stackrel\nu\lra& \End_R(A)\\
a\otimes b&\mapsto&(c\mapsto acb)
\end{eqnarray*}
is an isomorphism. Again, the Morita equivalence classes of 
Azumaya $R$-algebras form a group under tensor product over $R$, the
Azumaya Brauer group over $R$. Azumaya algebras are closely related to separable 
extension (cf \cite{AuslanderGoldman}).

\subsubsection{The commutative and graded case:}\label{commutativeandgradeAzumayasect}
If $R$ is a commutative $\Z$-graded ring, we say that a $\Z$-graded $R$-algebra $A$ 
is graded-Azumaya if $A$ is a finitely generated projective generator in the 
category of graded $R$-modules, and 
\begin{eqnarray*}
A\otimes_RA^{op}&\lra& \End_R(A)\\
a\otimes b&\mapsto&(c\mapsto acb)
\end{eqnarray*}
is a graded isomorphism. Two graded-Azumaya algebras $A$ and $B$ are equivalent if there is a 
faithfully projective 
graded $R$-module $C_A$ and faithfully projective 
graded $R$-module $C_B$ such that 
$$A\otimes_R\End_{R,graded}(C_A)\simeq B\otimes_R\End_{R,graded}(C_B).$$
Again, the equivalence classes form a group, the graded Brauer group $\ogrBr(R)$
with group law the tensor product $-\otimes_R-$ over $R$.

\subsubsection{The graded-commutative graded case:}
\label{gradedcommutativegradedSection}
For a graded commutative $\Z$-graded ring $R$ one considers $\Z$-graded $R$-algebras $A$. We call these 
rings $\Z$-graded-Azumaya (or graded-Azumaya for short)
if $A$ is a finitely generated projective generator in the category 
of graded $R$-modules, and 
\begin{eqnarray*}
A\widehat\otimes_RA^{op}&\lra& \End_R(A)\\
a\otimes b&\mapsto&(c\mapsto (-1)^{|b||c|}acb)
\end{eqnarray*}
is a graded isomorphism (recall the Koszul sign rule).
Again, two graded-Azumaya algebras $A$ and $B$ are equivalent if there is a 
faithfully projective 
graded $R$-module $C_A$ and faithfully projective 
graded $R$-module $C_B$ such that 
$$A\widehat\otimes_R\widehat\End_{R,graded}(C_A)\simeq B\widehat\otimes_R\widehat\End_{R,graded}(C_B).$$
Again, the equivalence classes form a group, the graded Brauer group $\grBr(R)$
with group law the tensor product $\widehat\otimes_R$ over $R$. 
The link between $\grBr(R)$ and $\ogrBr(R)$ is the same as the link between the 
Brauer group and the Brauer-Wall group of an ordinary field. 
If $K$ is concentrated in even degrees only, then a structure result is known.

By \cite[(7.8) Theorem]{CharlesSmall}, if $A$ is a $\Z/2\Z$-graded-Azuyama algebra
over an ungraded commutative ring $K$, then there is a short exact sequence
$$\xymatrix{
0\ar[r]&\textup{Br}(K)\ar[r]&\textup{BrW}(K)\ar[r]&Q_2(K)\ar[r]&0
}$$
where $Q_2(K)$ is an abelian group of exponent at most $4$, the group of quadratic extensions of $K$. 
For a $\Z$-graded ring $K$, concentrated in even degrees only, forgetting the grading, 
we let $\upsilon K$ be the ungraded ring $K$.
Likewise, Tilborghs and van Oystaeyen \cite{TilborghsVanOystaeyen}, generalized by 
Tilborghs in \cite{Tilborghs}, 
show that there is a commutative diagram with rows being short exact sequences 
$$
\xymatrix{
0\ar[r]&\textup{Br}(\upsilon K)\ar[r]&\textup{BrW}(\upsilon K)\ar[r]&Q_2(\upsilon K)\ar[r]&0\\
0\ar[r]&\grBr^{even}(K)\ar[r]\ar[u]^{\phi^{even}}&\grBr(K)\ar[u]^\phi\ar[r]&
Q_2^\Z(K)\ar@{^{(}->}[u]\ar[r]&0
}
$$
and where $\grBr^{even}(K)\simeq\ogrBr(K^{(1/2)})$. Here $\grBr^{even}(K)$
is the subgroup of $\grBr(K)$ generated by classes represented by 
graded-Azumaya algebras concentrated in even degrees only. Further,
$\ogrBr(K^{(1/2)})$ is the $\Z$-graded-Brauer group of Caenepeel and van-Oystaeyen
over the truly commutative graded-field $K^{(1/2)}$ which is $K$ but with a
different grading. Namely $(K^{(1/2)})_m=K_{2m}$. 
The group $Q_2^\Z(K)$ is the subgroup of $\Z$-graded quadratic extensions.\footnote{We actually proved these statements independently, before we discovered that they 
were proved by Tilborghs and van Oystaeyen \cite{TilborghsVanOystaeyen} a long 
time ago already.}
For more ample details we refer to these
papers. The right most maps in the short exact sequences are given by centralizer algebras,
basically.

\begin{Rem}
In \cite{dgBrauer} the second author defined a dg-Brauer group over a field $K$ 
as equivalence classes of central simple $K$-algebras which happen to be 
differential graded algebras. Equivalence was given analogous to the case of 
graded Brauer groups. Our $\dgBr^{II}(K,d_K)$ generalizes this case.
\end{Rem}

\subsection{Dg-division algebras}

A dg-module $(M,\delta)$ is dg-simple if it does not contain a dg-submodule other than $0$ and itself.

\begin{Def}
A dg-algebra $(A,d)$ is 
\begin{itemize}
\item a {\em dg-division algebra} if it does neither contain a left dg-ideal and nor a right dg-ideal 
other than $0$ or $A$.
\item a {\em dg-field} if it is commutative or graded commutative dg-division algebra. 
\end{itemize}
\end{Def}

We shall use frequently the following basic result.

\begin{Theorem} \cite{dgfields,dgseparable}\label{dgfieldcharacterisation}
Let $(A,d)$ be a dg-algebra. 
Then 
\begin{itemize}
\item
$(A,d)$ is a dg-division algebra if and only if $\ker(d)$ is a $\Z$-gr-division algebra (cf \cite{gradedrings}).
\item 
If $(A,d)$ is a dg-division algebra,
\begin{itemize} 
\item then,
\begin{itemize}
\item  either $d=0$ and $\ker(d)$ is  a graded-division algebra,
\item or $H(A,d)=0$ and there is a skew field $R_0$
such that $\ker(d)\simeq R_0[X,X^{-1};\phi]$ for an automorphism $\phi$ of $R_0$ 
and $Xr=\phi(r)X$ for any $r\in R_0$. 
\end{itemize}
\item If $H(A,d)=0$, then there is a homogeneous element $y\in A$ with
$d(y)=1$ and  $z=y^2\in\ker(d)$, and there is a map $D:\ker(d)\lra\ker(d)$  of degree $1$
defined by  $$D(a)=-(-1)^{|a|}d(yay)=ya-(-1)^{|a|}ay$$ 
for any homogeneous $a\in\ker(d)$,  
such that  $A$ is isomorphic with the quotient of the twisted polynomial ring
$$A\simeq\ker(d)[T;D]/ (T^2-z).$$ 
Moreover, the algebra structure on the twisted group ring is 
given by $D(a)=Ta-(-1)^{|a|}aT$ for any homogeneous $a\in\ker(d)$. Furthermore, 
$A=\ker(d)\oplus y\ker(d)$, and the isomorphism 
is 
\begin{eqnarray*}
\Phi:\ker(d)[T;D]&\lra& A\\
b+Ta&\mapsto&b+ya
\end{eqnarray*} 
for any homogeneous $a,b\in\ker(d)$. Further, 
for any homogeneous $a,b\in\ker(d)$ we get $d(b+ya)=a$. 
\end{itemize}
\end{itemize}
\end{Theorem}

\section{Classsification of commutative and graded-commutative dg-fields}
\label{dgfieldsclassification}

As a complement to Theorem~\ref{dgfieldcharacterisation},
as our first main result 
we give an explicit classification of dg-fields.

\begin{Theorem}\label{fieldclassscases}
Let $K$ be a commutative ring and let $(A,d)$ be a dg-field (commutative or graded commutative).
Then precisely one of the situations occur, and all these provide dg-fields.
\begin{enumerate}
\item \label{fieldclassscases1} $A$ is concentrated in degree $0$ and is a field.
\item \label{fieldclassscases2} $d=0$ and $A=L[T,T^{-1}]$ is a graded-field, 
where $T$ is in non zero degree and $L$ is a field (if $A$ is graded-commutative 
and the characteristic is different from $2$, then 
$|T|\in 2\Z$).
\item \label{fieldclassscases3} $\ker(d)_0=:L$ is a field, $A=L[T,T^{-1}]$ 
and $d(T^{2n+1})=T^{2n}$, $d(T^{2n})=0$ for all $n$.
\item \label{fieldclassscases4} $\ker(d)_0=:L$ is a field, 
$A=\ker(d)\otimes_KK[Y]/Y^2$, $d(Y)=1$, and 
    \begin{enumerate}
    \item either $\ker(d)=L$ 
    \item or $\ker(d)=L[T,T^{-1}]$ for some $T$ in non zero even degree.
    \end{enumerate} 
\item \label{fieldclassscases5}  $\ker(d)_0=:L$ is a field, 
$A=L[T,T^{-1},Y]/(Y^2-\lambda T^2)$ is commutative, for some $\lambda\in L$, 
$d(Y)=1$, $d(T^n)=0$ for all $n$, $T$ and $Y$ are of degree $-1$. Note that $\lambda$ can be $0$ or 
can be different from $0$.
\end{enumerate}
\end{Theorem}

Proof. We know that $\ker(d)$ is a graded-field, and hence 
either $\ker(d)$ is concentrated in degree $0$, or else 
$\ker(d)=L[T,T^{-1}]$ for some field $L$, and some variable $T$ in non zero degree.
The case $d=0$ then provides the cases (\ref{fieldclassscases1}) and 
(\ref{fieldclassscases2}) in the statement.

\medskip

In the sequel we may suppose that $d\neq 0$, and hence $(A,d)$ is acyclic, 
which implies directly that $A_{-1}\neq 0$.

\medskip

{\em Suppose first that $\ker(d)$ is concentrated in even degrees only. }
We have two cases.

(i) In the first case $\ker(d)=\ker(d)_0$ is a field, 

(ii) and in the second case  
$\ker(d)=\ker(d)_0[T,T^{-1}]$ for some $T$ of degree $n\in 2\Z$.

Consider a homogeneous element $y'\in A$ of degree $-1$. It cannot be in 
$\ker(d)$, since $\ker(d)$ is concentrated in even degrees only. 
Whence $d(y')=u\in A_0\setminus\{0\}$, and actually $u\in\ker(d)_0$.
Hence, $u$ is invertible, and for $y:=u^{-1}y'$ we get $d(y)=1$.
If there is another $\tilde y$ in degree $-1$,  
with $y$, then $d(\widetilde y)=v\neq 0$ and 
$vy-\widetilde y\in\ker(d)$. Since $\ker(d)$ does not allow any element of odd 
degree, $A_{-1}=\ker(d)_0\cdot y$. 

Suppose $y^2=0$. 
Now, in the first case (i) we get 
$A=\ker(d)_0[Y]/Y^2=\ker(d)[Y]/Y^2$ with $d(Y)=1$, by mapping $Y$ to $y$. 

In the second case (ii) we have that $T$ is invertible. Let $b\in A_m\setminus\ker(d)$ be 
an element of degree $m$. Then $d(b)\in\ker(d)$, and hence $d(b)=uT^n$ for some $n$
and some $0\neq u\in\ker(d)_0$. Then, $T^{-n}b\in A_{-1}=\ker(d)_0\cdot y$ and $m=n\cdot|T|-1$. 
Therefore, again mapping $Y$ to $y$,
$$A=\ker(d)_0[Y]/Y^2\otimes_{\ker(d)_0}\ker(d)=\ker(d)[Y]/Y^2$$
in this case as well, and hence in both cases. This gives the 
item (\ref{fieldclassscases4}) in the statement of the theorem.

Now suppose that $y^2\neq 0$.
Note that $z=y^2\in\ker(d)$. 
Then, by Aldrich-Garcia Rozas' theorem we get, mapping $Y$ to $y$ again, that 
$$A={\ker(d)[Y]/(Y^2-z)},$$
taking into account that $A$ was supposed to be commutative or graded commutative, 
and hence the derivation $D$ in the theorem is $0$. 
Since $z\neq 0$, we get that $z$ is invertible in $\ker(d)$. Hence,  
$$L[z,z^{-1}]=\ker(d)\textup{ and }A=\ker(d)[Y]/(Y^2-z)=L[z,z^{-1},Y]/(Y^2-z).$$ 
This then yields the situation (\ref{fieldclassscases3}) of the statement 
of the theorem with $Y$ playing the role of $T$.

\medskip

{\em Suppose now that $\ker(d)$ has a non zero component in odd degrees.}
Since $\ker(d)$ is a graded-field, and for any $b\in\ker(d)$ of odd degrees, 
we would get that $b^2=0$ in case $(A,d)$ is graded-commutative 
(hence of characteristic different from $2$), $b$ cannot be invertible in 
this case. Hence, in this case we have that $(A,d)$ needs to be commutative actually.
Still $\ker(d)$ is a graded-field, and hence $\ker(d)=\ker(d)_0[T,T^{-1}]$ for some $T$ of 
degree $2m+1$. Since the algebra was assumed to be 
commutative, the derivation $D$ in  
Aldrich-Garcia Rozas' theorem is $0$ on even degree elements and 
$D(a)=2aT$ on odd degree elements, producing $Ta=aT$ for any $a\in \ker(d)$.
Hence, by Aldrich-Garcia Rozas' theorem
there is an element $y\in A_{-1}$ with $d(y)=1$, satisfying $z=y^2\in\ker(d)$, and 
mapping $Y$ to $y$ again,
$$A=\ker(d)[Y]/(Y^2-z)=\ker(d)_0[T,T^{-1},Y]/(Y^2-z)$$  
for some $Y$ of degree $-1$. Since $y^2\in\ker(d)$, either $y^2=0$, or else $z=y^2$ is in degree $-2$.
If $z\neq 0$, we have that $\ker(d)_{-2}\neq 0$, and therefore $T$ is of degree $-1$. 
Since $z\in\ker(d)_{-2}=\ker(d)_0\cdot T^2$, there is $\lambda\in\ker(d)_0$ such that
$z=\lambda T^2$. Note that $\lambda$ does not need to be a square in $\ker(d)_0$.
We get  the statement of case (\ref{fieldclassscases5}) of the theorem. 

Altogether this then proves the theorem. \dickebox

\begin{Rem}
For a dg-field $(K,d_K)$ we have that $\ker(d_K)$ is a graded-field 
(cf \cite{dgfields} and \cite{dgseparable}).
Any graded-module over a graded-field allows a graded-basis (cf \cite[Lemma 1.7]{vanGeel}), and hence 
is graded-flat.  
\end{Rem}

\begin{Rem}
Note that the classification of 
Proposition~\ref{fieldclassscases} shows that any graded-commutative 
dg-division algebra is actually commutative.
\end{Rem}

\section{Defining dg-Brauer groups}
\label{definingdgbrauersection}

\subsection{Definition of dg-Azumaya algebras}
Let $K$ be a commutative ring. 
Recall that an Azumaya algebra over  $K$ is a $K$-algebra $A$, 
which is faithfully projective over $K=Z(A)$, 
and such that the map 
\begin{eqnarray*}
A\otimes_KA^{op}&\stackrel\nu\lra&\End_K(A)\\
a\otimes b&\mapsto&(x\mapsto axb)
\end{eqnarray*}
is an isomorphism.

Similarly a graded-Azuyama algebra is a graded algebra $A$ which is 
faithfully projective over $K=Z_{gr}(A)$  and such that 
the map 
\begin{eqnarray*}
A\widehat\otimes_KA^{op}&\stackrel\mu\lra&\widehat\End_{K,graded}(A)\\
a\otimes b&\mapsto&(x\mapsto (-1)^{|b||x|}axb)
\end{eqnarray*}
is an isomorphism. Note that we use Koszul signs in order to be able to cope with the 
various signs appearing in the construction 
of dg-algebras. 
In particular, we shall need to assume that the base ring $R$ is graded-commutative.
By \cite[Proposition III.4.1]{caenepeelvanoystaeyen} 
we see that a graded algebra $A$ over a commutative $R$ 
is graded-Azumaya if and only if $A$ is Azumaya.

\begin{Def}

\begin{itemize}
\item Let $(K,d_K)$ be a graded commutative dg-ring and let $(A,d)$ be a dg-algebra.
We say that $(A,d)$ is {\em dg-Azumaya of the first kind} if $\ker(d)$ is graded-Azumaya
algebra over $\ker(d_K)$.

We say that two dg-Azumaya algebras $(A,d_A)$ and $(B,d_B)$ of the first kind 
are equivalent if $\ker(d_A)$ and $\ker(d_B)$ represent the same class 
in the graded-Brauer group $\grBr(\ker(d_K))$.
\item
Let $(K,d_K)$ be a graded commutative dg-ring and let $(A,d)$ be a dg-algebra.
We say that $(A,d)$ is {\em dg-Azumaya of the second kind} if $A$ is graded-Azumaya
algebra over $K$.

We say that two dg-Azumaya algebras $(A,d_A)$ and $(B,d_B)$ of the second kind 
are equivalent if there are dg-modules $(C_A,\partial_A)$ and $(C_B,\partial_B)$, which are 
faithfully projective as graded $K$-modules, and such that 
$$(A,d_A)\widehat\otimes_K (\End_K^\bullet(C_A),d_{\Hom})\simeq 
(B,d_B)\widehat\otimes_K(\End_K^\bullet(C_B),d_{\Hom})$$
as dg-algebras.
\end{itemize}
\end{Def}

\subsection{The group structure of dg-Brauer groups of the first and second kind}

\begin{Prop}\label{groupstructureofdgbrauerofsecondkind}
Let $(K,d_K)$ be a graded commutative dg-ring. Then the set of equivalence classes of 
dg-Azumaya algebras of the second kind over $(K,d_K)$ form a group 
$\dgBr^{II}(K,d)$ with group law being given by 
$$[(A,d_A)]\cdot[(B,d_B)]:=[(A\widehat\otimes_KB,d_{A\widehat\otimes_KB})].$$
\end{Prop}

Proof. 
The group law is clearly well-defined. 
We have that the map 
\begin{eqnarray*}
A\widehat\otimes_KA^{op}&\stackrel\mu\lra&\widehat\End_{K,graded}(A)\\
(a\otimes b)&\mapsto&(x\mapsto(-1)^{|b||x|}axb)
\end{eqnarray*}
is an isomorphism. We need to check that the isomorphism maps the differential on the tensor 
product to $d_{\Hom}$ on the complex $(A,d_A)$.
\begin{eqnarray*}
d_{\Hom}(\mu(a\otimes b))(x)
&=&(d_A\circ\mu(a\otimes b)-(-1)^{|a|+|b|}\mu(a\otimes b)\circ d_A)(x)\\
&=&d_A((-1)^{|b||x|}axb)-(-1)^{|a|+|b|+(|b|(|x|+1))}ad_A(x)b\\
&=&(-1)^{|b||x|}\left(d_A(axb)-(-1)^{|a|}ad_A(x)b\right)\\
&=&(-1)^{|b||x|}\left(d_A(a)xb+(-1)^{|a|}ad_A(xb)-(-1)^{|a|}ad_A(x)b\right)\\
&=&(-1)^{|b||x|}\left(d_A(a)xb+(-1)^{|a|}ad_A(x)b+(-1)^{|a|+|x|}axd_A(b)-(-1)^{|a|}ad_A(x)b\right)\\
&=&(-1)^{|b||x|}\left(d_A(a)xb+(-1)^{|a|+|x|}axd_A(b)\right)\\
&=&(-1)^{|b||x|}d_A(a)xb+(-1)^{|a|+(|b|+1)|x|}axd_A(b)\\
&=&\mu(d_A(a)\otimes b+(-1)^{|a|}a\otimes d_A(b))(x)\\
&=&\mu((d_A\otimes 1+1\otimes d_A)(a\otimes b))(x)
\end{eqnarray*}
Hence,
$$[(A,d_A)]\cdot[(A^{op},d_A)]=[(\End_K^\bullet(A),d_{\Hom})]=[(K,d_K)]$$
by the definition of a graded-Azumaya algebra and the equivalence relation. 
Associativity of the tensor product is clear.
This shows the lemma. \dickebox

\begin{Lemma} \label{tensoroftwodegazuyamaoffirstkind}
Let $(K,d_K)$ be a graded-commutative dg-algebra.
Let $(A,d_A)$ and $(B,d_B)$ be dg-Azumaya algebras of the first kind over $(K,d_K)$. 
Then $$(\ker(d_A)\widehat\otimes_{\ker(d_K)}\ker(d_B))\widehat\otimes_{\ker(d_K)}K$$
is a dg-Azumaya algebra of the first kind. 
\end{Lemma}

Proof. By the arguments in \cite[III.4.5]{caenepeelvanoystaeyen}, replacing 
the ordinary tensor product by the graded-tensor product $\widehat\otimes$ 
and the ordinary $\textup{End}$ by the graded $\widehat\End_{graded}$, we see that 
$(\ker(d_A)\widehat\otimes_{\ker(d_K)}\ker(d_B))$ is a graded-Azumaya $\ker(d_K)$-algebra.
By definition, a graded-Azumaya algebra over $\ker(d_K)$ is 
faithfully projective as a graded  $\ker(d_K)$-module.
Hence,  $(\ker(d_A)\widehat\otimes_{\ker(d_K)}\ker(d_B))=:X$
is projective, and therefore flat as a  $\ker(d_K)$-module. This implies that
$$0\lra X\widehat\otimes_{\ker(d_K)}\ker(d_K)\lra 
X\widehat\otimes_{\ker(d_K)}K\lra X\widehat\otimes_{\ker(d_K)}d_K(K)[1]\lra 0$$
is exact.
\begin{eqnarray*}
\ker(\id_{\ker(d_A)\widehat\otimes_{\ker(d_K)}\ker(d_B)}\widehat\otimes_{\ker(d_K)}d_K)
&=&\ker(d_A)\widehat\otimes_{\ker(d_K)}\ker(d_B)\widehat\otimes_{\ker(d_K)}\ker(d_K)\\
&=&\ker(d_A)\widehat\otimes_{\ker(d_K)}\ker(d_B)
\end{eqnarray*}
using that $X$ is graded-flat as a $\ker(d_K)$ module.
This shows the lemma. \dickebox

\begin{Prop}\label{grouplawindgbrI}
Let $(K,d_K)$ be a graded-commutative dg-algebra. Then the
equivalence classes of dg-Azumaya of the first kind form a group $\dgBr^{I}(K,d_K)$
by the group law 
$$[(A,d_A)]\cdot [(B,d_B)]:=[(\ker(d_A)\widehat\otimes_{\ker(d_K)}\ker(d_B)
\widehat\otimes_{\ker(d_K)}K,\id\widehat\otimes d_K)]$$
\end{Prop}

Proof. Lemma~\ref{tensoroftwodegazuyamaoffirstkind} shows that 
$\ker(d_A)\widehat\otimes_{\ker(d_K)}\ker(d_B)\widehat\otimes_{\ker(d_K)}K$ is a 
dg-Azumaya algebra of the first kind. 

We shall first verify that the group law is well-defined. 
If $[(A,d_A)]=[(A',d_{A'})]$, then, by definition, 
$[\ker(d_A)]_{graded}=[\ker(d_{A'})]_{graded}.$ But as is shown in  
\cite[III.4.5]{caenepeelvanoystaeyen}, trivially adapted to the 
graded situation, the group law in the graded Brauer group 
is given by $-\widehat\otimes_{\ker(d_K)}-$, and this is well-defined with respect to 
equivalence of graded algebras.  Similarly, we argue for the second component. 

Further, $[(A^{op},d_A)]$ is the inverse to $[(A,d_A)]$, since 
the cycles of   $[(A^{op},d_A)]$ are precisely $\ker(d_A)^{op}$, and its 
graded equivalence class is the inverse of $[\ker(d_A)]_{graded}$.
Associativity is clear, and we showed the result. \dickebox 

\section{Computing dgBrauer groups of the first and second kind}

\label{differentialzerosection}

\subsection{dgBrauer group of the second kind and graded Brauer groups, 
the differential zero case} 

\begin{Rem}\label{goodgrading}
Commutative or graded commutative graded-fields $R$ behave very much like ordinary fields. 
All graded modules over such rings are free and allow bases (cf \cite[Lemma 1.7]{vanGeel}).  
Further, by \cite{BeattieLiuHong} a Chevalley-Jacobson density theorem holds. These are the
only ingredients in the proof of \cite[Theorem 1.4]{DascalescuIonNastasescuRioMontes}.
Hence, the proof of this result can be copied literally to see that any
matrix algebra over a commutative or graded commutative graded-fields $R$
is graded-isomorphic to one with a good grading, in the sense that 
the elementary matrices are homogeneous. This in turn is the hypothesis for
our proof from \cite[Proposition 3.3]{dgBrauer} that any differential 
on $\End_{graded,K}(P)$ is 
actually equivalent to $d_{\Hom}$ for some differential $\partial$ on $P$, 
where $K$ is a field.
This last statement is therefore true as well for commutative or 
graded commutative graded-fields $R$,
with literally the same proof as for ungraded fields. 
\end{Rem}

\begin{Theorem}\label{dgBrIIvsusgradedBrK}
Let $(K,d_K)$ be a graded commutative dg-ring. Then there is a group homomorphism 
$$\Psi:\dgBr^{II}(K,d_K)\lra \grBr(K).$$

If $d_K=0$, then  $$\dgBr^{II}(K,0)\simeq L\times \grBr(K)$$
for the subgroup $L$ of $\dgBr^{II}(K,0)$ formed by the class of dg-algebras obtained by
putting a differential on the graded endomorphism algebra 
of some faithfully graded projective $K$-module. 

If $K$ is a graded-field, then $L=0$.
\end{Theorem}

Proof.
Let $[(A,d_A)]$ be an equivalence class in $\dgBr^{II}(K,d_K)$. Then $A$ is graded-Azumaya, by definition. 
Further, if $[(A,d_A)]=[(B,d_B)]$ in $\dgBr^{II}(K,d_K)$, then there are dg-modules  
$(C_A,\partial_A)$ and $(C_B,\partial_B)$ over $(K,d_K)$, faithfully projective over $K$, 
such that 
$$(A,d_A)\widehat\otimes_K\End_K^\bullet(C_A,\partial_A)\simeq (B,d_B)\widehat\otimes_K\End_K^\bullet(C_B,\partial_B)$$
as dg-algebras.
Forgetting the differential this implies 
$$A\widehat\otimes_K\widehat\End_{K,graded}(C_A)\simeq B\widehat\otimes_K\widehat\End_{K,graded}(C_B)$$
as graded algebras,
and hence $[A]=[B]$ in $\grBr(K)$. Therefore, forgetting the differential
yields a group homomorphism $$\Psi:\dgBr^{II}(K,d_K)\lra \grBr(K).$$

Suppose now $d_K=0$.
Then $\Psi$ is a split epimorphism. 
Indeed, let $A$ be a graded-Azumaya algebra over $K$. Then  $(A,0)$ is a dg-algebra, since $K$ 
is equipped with differential $0$. Moreover, $(A,0)$ is dg-Azumaya 
of the second kind over $K$. Further, if $[A]=[B]$, then 
there are faithfully projective graded modules  
$C_A$ and $C_B$ over $K$, such that 
$$A\widehat\otimes_K\widehat\End_{K,graded}(C_A)\simeq B\widehat\otimes_K\widehat\End_{K,graded}(C_B)$$
as graded algebras. Since the differential on $K$ is $0$, the graded modules $(C_A,0)$ and $(C_B,0)$ are
dg-modules over $(K,0)$. 
Hence, we get an isomorphism 
$$(A,0)\widehat\otimes_K\End_K^\bullet(C_A,0)\simeq (B,0)\widehat\otimes_K\End_K^\bullet(C_B,0)$$
of dg-algebras. Hence, 
\begin{eqnarray*}
\grBr(K)&\stackrel{\Phi}\lra&\dgBr^{II}(K,0)\\
\ [A]&\mapsto&[(A,0)]
\end{eqnarray*}
is a group homomorphism satisfying $\Psi\circ\Phi=\id$.
In particular, $\Psi$ is surjective. 
 
Let $(A,d)$ be a dg-algebra such that $A$ is graded-Azumaya over $K$
in the kernel of $\Psi$. Hence, $A$ as a graded algebra is equivalent to 
a graded endomorphism algebra of a faithful graded projective graded $K$-module $P$.
The differential $d$ is then a differential on the graded endomorphism algebra 
of a faithfully projective $K$-module.
But this precisely means that $(A,d)$ belongs to the subgroup defined by $L$.
Hence 
$\dgBr^{II}(K,0)\simeq L\sdp \grBr(K)$ for $L=\ker(\Psi)$. However, 
$\dgBr^{II}(K,0)$ is an abelian group (cf \cite[Lemma 2.6]{dgBrauer}). 
Therefore the semidirect product 
is a direct product and 
$$\dgBr^{II}(K,0)\simeq L\times \grBr(K)$$ for $L=\ker(\Phi)$.

Now, suppose that $K$ is a graded-field. Then by 
\cite[Proposition IV.1.3]{caenepeelvanoystaeyen}
we get that $A$ is graded-simple central.
By Remark~\ref{goodgrading} we obtain that $(A,d)$ is equivalent to $(K,0)$. 

This shows the theorem. 
\dickebox

\begin{Rem}
It would be nice to get a necessary and sufficient criterion when $L=0$.
\end{Rem}

\subsection{dgBrauer groups of the first kind and graded Brauer groups}

dgBrauer groups of the first kind are designed to be actually isomorphic to the 
graded Brauer group of the cycles, as is observed by the following easy lemma.

\begin{Lemma}\label{dgBrIisgrBr}
If $(K,d_K)$ is a  graded commutative  dg-algebra, 
then $\dgBr^I(K,d_K)\simeq \grBr(\ker(d_K))$. 
\end{Lemma}

Proof. 
We shall use 
Proposition~\ref{grouplawindgbrI}. 

Let $(A,d_A)$ be a dg-Azumaya algebra of the first kind. Then 
by definition $\ker(d_A)$ is a graded-Azuyama algebra over $\ker(d_K)$. 
Further, by definition, equivalent dg-Azumaya algebras of the first 
kind yield equivalent graded algebras of cycles. 
Hence, taking cycles yields a well-defined group homomorphism 
$$\dgBr^I(K,d_K)\lra \grBr(\ker(d)).$$
This map is injective since two dg-Azumaya algebras of the first kind 
with the same image in $\grBr(\ker(d_K))$ yield 
equivalent classes in $\dgBr^I(K,d_K)$, by definition. 

Given a graded-Azumaya algebra $C$ over $\ker(d_K)$, then $C\widehat\otimes_{\ker(d_K)}K$
is a dg-algebra with differential $\partial_C:=\id_C\widehat\otimes d_K$. Since  $C$ is 
graded Azumaya over $\ker(d_K)$, it is finitely generated graded-projective 
over $\ker(d_K)$, and hence flat over $\ker(d_K)$.
This implies that $\ker(\partial_C)=C$. Hence, the above group homomorphism is surjective 
as well. This proves the lemma. \dickebox

\begin{Rem}
Suppose that $(K,d_K)$ is a commutative dg-ring. 
Then, for two dg-$R$-algebras $(A,d_A)$ and $(B,d_B)$ such that $\ker(d_A)$ and  
$\ker(d_B)$ are ordinary-graded algebras over $\ker(d_K)$, then
we may consider the 
ordinary tensor product on cycles $\ker(d_A)\otimes_{\ker(d_K)}\ker(d_B)$ 
and define a dg-Brauer group of the first kind this way. All
statements in this subsection on $\dgBr^I(K)$ remain true, after appropriate modifications, 
and we get an isomorphism of the so-defined dg-Brauer group with $\ogrBr(\ker(d_K))$. 
\end{Rem}

\subsection{Linking dgBrauer groups of the first and the second kind; the acyclic case}

In this section we shall consider the case of $(K,d_K)$ being acyclic. Then,
by the main theorem in \cite{Tempest-Garcia-Rochas} we have that every 
dg-algebra over $(K,d_K)$ is acyclic as well. We shall use the following 
result frequently.

\begin{Lemma}\label{acyclicalgebrasareinducedfromcycles}
Let $(K,d_K)$ be an acyclic  graded commutative dg-algebra, 
and let $(A,d_A)$ be a dg-algebra over $(K,d_K)$. Then there is an 
isomorphism of dg-algebras
$$
 \ker(d_A)\widehat\otimes_{\ker(d_K)}(K,d_K)\stackrel{\alpha}{\lra} (A,d_A)
$$
given by $a\otimes x\mapsto ax$.
\end{Lemma}

Proof. 
We get 
\begin{eqnarray*}
\alpha((a\otimes x)\cdot(b\otimes y))&=&\alpha((-1)^{|x||b|}ab\otimes xy)\\
&=&(-1)^{|x||b|}abxy\\
&=&axby\\
&=&\alpha(a\otimes x)\cdot\alpha(b\otimes y)
\end{eqnarray*}
Hence, $\alpha$ is an algebra homomorphism. 
Further, denote by $\partial_{\ker(d_A)\otimes K}=\id\otimes d_K$ the differential on $\ker(d_A)\otimes_{\ker(d_K)}K$,
\begin{eqnarray*}
\alpha(\partial_{\ker(d_A)\otimes K}(a\otimes x))&=&
\alpha((-1)^{|a|}a\otimes d_K(x))\\
&=&(-1)^{|a|}a\cdot d_K(x)\\
&=&d_A(a\cdot x)\\
&=&d_A(\alpha(a\otimes x))
\end{eqnarray*}
using that $a\in\ker(d_A)$.

Since $(K,d_K)$ is acyclic, also $(A,d_A)$ is acyclic, By \cite{Tempest-Garcia-Rochas}
we get that 
$$A=\ker(d_A)\oplus\ker(d_A)T_A=\ker(d_A)\oplus T_A\ker(d_A)$$ 
as left (resp. right) modules over $\ker(d_A)$, and also 
$$K=\ker(d_K)\oplus T_K\ker(d_K)=\ker(d_K)\oplus \ker(d_K)T_K$$
for some $T_K\in K$ and $d_K(T_K)=1$. Observe that any homogeneous $T_A$ with $d(T_A)=1_A$ 
does have this property. 
As $(A,d_A)$ is a dg-module over $(K,d_K)$, we obtain by definition a dg-ring homomorphism 
$$(K,d_K)\stackrel{\nu}\lra Z_{gr}(A,d_A)$$
since $K$ is graded commutative.
This implies $$d_A(\nu(T_K))=\nu(d_K(T_K))=\nu(1)=1_A$$
and hence we may identify $$\nu(T_K)=T_A\in Z_{gr}(A).$$
Hence, we pose $T:=T_A=\nu(T_K)$, the image of $T_K$ in $A$. Therefore, 
$$A\simeq \ker(d_A)\oplus\ker(d_A)T\simeq \ker(d_A)\otimes_{\ker(d_K)}K.$$
This shows that $\alpha$ is bijective and we proved the lemma. \dickebox

\begin{Prop} \label{secondkindimpliesfirstkind}
Let $(K,d_K)$ be an acyclic and graded commutative dg-ring.
Let $(A,d_A)$ be a dg-Azumaya algebra of the second kind over $(K,d_K)$. 
If $\ker(d_A)$ is faithfully projective as a $\ker(d_K)$-module, then 
$(A,d_A)$ is a dg-Azumaya algebra of the first kind. 
\end{Prop}

Proof. 
Suppose that the map 
\begin{eqnarray*}
A\widehat\otimes_KA^{op}&\stackrel\phi\lra&\widehat\End_{K,graded}(A)\\
a\otimes b&\mapsto&(x\mapsto (-1)^{|b||x|}axb)
\end{eqnarray*}
is an isomorphism.

Now, 
$$A\widehat\otimes_KA^{op}\simeq \widehat\End_{K,graded}(A)$$
induces a differential $\partial_E$ on $\widehat\End_{K,graded}(A)$
by transport of structure, from the differential $d_A\otimes d_{A^{op}}$ on $A\otimes A^{op}$.
Again, this way $(\widehat\End_{K,graded}(A),\partial_E)$ is a dg-module over $(K,d_K)$, and 
is hence acyclic as well. 
By Lemma~\ref{acyclicalgebrasareinducedfromcycles} and its proof,
\begin{eqnarray*}
A\widehat\otimes_KA^{op}&\simeq&\left(\ker(d_A)\widehat\otimes_{\ker(d_K)}K\right)
\widehat\otimes_K
\left(K\widehat\otimes_{\ker(d_K)}\ker(d_A)\right)^{op}\\
&\simeq&\ker(d_A)\widehat\otimes_{\ker(d_K)}K\widehat\otimes_{\ker(d_K)}\ker(d_A)^{op}\\
&\simeq&\left(\ker(d_A)\widehat\otimes_{\ker(d_K)}\ker(d_A)^{op}\right)\oplus 
\left(\ker(d_A)\widehat\otimes_{\ker(d_K)}\ker(d_A)^{op} \cdot T\right)\\
&\stackrel\chi\lra&\widehat\End_{\ker(d_K),graded}(\ker(d_A))\oplus \widehat\End_{\ker(d_K),graded}(\ker(d_A))T
\end{eqnarray*}
where $\chi$ is the homomorphism given by left and right multiplication. 
We know by hypothesis that $\chi$ is injective, and surjective to 
$\widehat\End_{K,graded}(A)$.
However, 
\begin{eqnarray*}
\widehat\End_{K,graded}(A)&=&
\widehat\Hom_{\ker(d_K)\widehat\otimes_{\ker(d_K)}K,graded}
(\ker(d_A)\widehat\otimes_{\ker(d_K)}K,\ker(d_A)\widehat\otimes_{\ker(d_K)}K)\\
&=&\widehat\Hom_{\ker(d_K),graded}(\ker(d_A),\ker(d_A)\widehat\otimes_{\ker(d_K)}K)\\
&=&\widehat\Hom_{\ker(d_K),graded}(\ker(d_A),\ker(d_A))\oplus \widehat\Hom_{\ker(d_K),graded}(\ker(d_A),\ker(d_A))T
\end{eqnarray*}
and hence $\chi$ is an isomorphism. Now, the only non trivial action 
of the differential $\partial_E$ is on $T$, and there we get $\partial_E(T)=1$. 
But then, the isomorphism 
\begin{eqnarray*}
\left(\ker(d_A)\widehat\otimes_{\ker(d_K)}\ker(d_A)^{op}\right)\lefteqn{\oplus 
\left(\ker(d_A)\widehat\otimes_{\ker(d_K)}\ker(d_A)^{op} \cdot T\right)}\\
&\stackrel\phi\lra&\widehat\End_{\ker(d_K),graded}
(\ker(d_A))\oplus \widehat\End_{\ker(d_K),graded}(\ker(d_A))T
\end{eqnarray*}
yields an isomorphism of cycles, which gives 
$$\ker(d_A)\widehat\otimes_{\ker(d_K)}\ker(d_A)^{op}\simeq \widehat\End_{\ker(d_K),graded}(\ker(d_A)).$$
The fact that $\ker(d_A)$ is projective as a $\ker(d_K)$-module 
holds by assumption.
This shows the proposition. \dickebox

\medskip

By Aldrich and Garcia-Rozas' theorem we get for an acyclic 
dg-algebra $(A,d_A)$ an equivalence of categories 
$$F_A:(A,d_A)-dgmod\lra \ker(d_A)-grmod$$ 
given by taking cycles and a quasi-inverse $G_A$ given by $A\widehat\otimes_{\ker(d_A)}-$. 

\begin{Lemma}\label{Fcommuteswithtensor}
Let  $(K,d_K)$ be an acyclic graded commutative 
dg-ring, and let $(A,d_A)$ and $(B,d_B)$ be
dg-algebras over $(K,d_K)$. 
Then $$F_{A\widehat\otimes_KB}(A\widehat\otimes_KB)\simeq F_A(A)\widehat\otimes_{\ker(d_K)}F_B(B).$$
\end{Lemma}

Proof.
By Aldrich and Garcia-Rozas' theorem we get that $(A,d_A)$ and $(B,d_B)$ are both 
acyclic, and $$A\simeq \ker(d_A)\oplus \ker(d_A)T_A$$ as well as 
$$B\simeq \ker(d_B)\oplus \ker(d_B)T_B$$
as modules over $\ker(d_A)$, resp. $\ker(d_B)$ 
for $d_A(T_A)=1_A$ and $d_B(T_B)=1_B$ for homogeneous elements $T_A,T_B$. 
Any such element has this property. 
Therefore we may take $T_A$ the image of $T_K$ in $A$, and $T_B$ the image of $T_K$
in $B$, and, slightly abusing the notation, abbreviate all these elements by $T$.  
\begin{eqnarray*}
A\widehat\otimes_KB&\simeq& (\ker(d_A)\oplus \ker(d_A)T)\widehat\otimes_K(\ker(d_B)\oplus \ker(d_B)T)\\
&\simeq&(\ker(d_A)\widehat\otimes_{\ker(d_K)}K)\widehat\otimes_K(K\widehat\otimes_{\ker(d_K)}\ker(d_B))\\
&\simeq&\ker(d_A)\widehat\otimes_{\ker(d_K)}K\widehat\otimes_{\ker(d_K)}\ker(d_B)\\
&\simeq&\ker(d_A)\widehat\otimes_{\ker(d_K)}(\ker(d_K)\oplus\ker(d_K)T)\widehat\otimes_{\ker(d_K)}\ker(d_B)\\
&\simeq&(\ker(d_A)\widehat\otimes_{\ker(d_K)}\ker(d_B))\oplus (\ker(d_A)\widehat\otimes_{\ker(d_K)}\ker(d_B)T)
\end{eqnarray*}
Now, as $A\widehat\otimes_KB$ is a dg-module over $(K,d_K)$ again, $A\widehat\otimes_KB$ is acyclic, and 
the differential on the subalgebra 
$(\ker(d_A)\widehat\otimes_{\ker(d_K)}\ker(d_B))$ is $0$, whereas $d_{A\widehat\otimes_KB}(T)=1$.
This shows the lemma. \dickebox

\begin{Theorem}\label{theultimateforacyclic}
Let $(K,d_K)$ be an acyclic  graded-commutative dg-ring. 
Assume that if a dg-module $(P,\partial_P)$ is faithfully 
projective as a graded $K$-module, then
also $\ker(\partial_P)$ is faithful projective as a graded-module
over $\ker(d_K)$.
Then there an isomorphism
$$\dgBr^{II}(K,d_K)\lra \dgBr^{I}(K,d_K)$$
given by the identity on representatives of equivalence classes of objects 
and hence 
$$\dgBr^{II}(K,d_K)\simeq  \grBr(\ker(d_K)).$$
\end{Theorem}

\begin{Rem}
Note that if $(K,d_K)$ is a dg-field, then 
$\ker(d_K)$ is a graded-field, and modules 
over graded-fields are graded-free.
The technical assumption in Theorem~\ref{theultimateforacyclic} 
is true more generally if 
$(K,d_K)-dgmod$ is a semisimple category.
Note that $(K,d_K)-dgmod$ is a semisimple category implies that $(K,d_K)$ is acyclic by 
Aldrich and Garcia-Rozas' theorem \cite{Tempest-Garcia-Rochas}. 
\end{Rem}

Proof.
Let $(A,d)$ be a dg-Azumaya algebra of the second kind, by
Proposition~\ref{secondkindimpliesfirstkind}
we get that $(A,d)$ is a dg-Azumaya algebra of the first kind. 

We need to see that the map 
\begin{eqnarray*}
\dgBr^{II}(K,d_K)&\lra& \dgBr^{I}(K,d_K)\\
\ [(A,d)]&\mapsto&[(A,d)]
\end{eqnarray*}
is well defined and is a group homomorphism.

Let $(P,\partial_P)$ be a dg-$(K,d_K)$-module, which is faithfully projective as a 
graded $K$-module. Then any dg-endomorphism $\varphi$ of 
$(P,\partial_P)$ induces a graded $\ker(d_K)$-linear endomorphism of $\ker(\partial_P)$.
Indeed, by Aldrich and Garcia-Rozas theorem \cite{Tempest-Garcia-Rochas}
there is an equivalence of categories
\begin{eqnarray*}
gr-mod(\ker(d_K))&\lra&dg-mod(K,d_K)\\ 
M&\mapsto&K\widehat\otimes_{\ker(d_K)}M
\end{eqnarray*}
with quasi-inverse 
\begin{eqnarray*}
dg-mod(K,d_K)&\lra&gr-mod(\ker(d_K))\\ 
(N,\partial_N)&\mapsto&\ker(\partial_N).
\end{eqnarray*}
Hence, the above construction gives an isomorphism of algebras
$$
\End^\bullet_{(K,d_K)}(P,\partial_P)\simeq 
\widehat\End_{\ker(d_K),graded}(\ker(\partial_P))
$$
Also, since $P$ is faithfully projective as a $K$-module, 
$\ker(\partial_P)$ is faithfully projective as a $\ker(d_K)$-module.

Now, if $[(A,d_A)]=[(B,d_B)]$ in $\dgBr^{II}(K,d_K)$, then there are
dg-modules $(C_A,\partial_A)$ and $(C_B,\partial_B)$, faithfully projective over $K$,
such that
$$(A,d_A)\widehat\otimes_K\End_K^\bullet((C_A,\partial_A),d_{\Hom}^A)\simeq 
(B,d_B)\widehat\otimes_K\End_K^\bullet((C_B,\partial_B),d_{\Hom}^B)$$
as dg-algebras. 
Using now Lemma~\ref{Fcommuteswithtensor}, we get that, taking 
cycles yields an isomorphism of graded algebras
$$\ker(d_A)\widehat\otimes\ker(d_{\Hom}^A)\simeq \ker(d_B)\widehat\otimes\ker(d_{\Hom}^B).$$
Now, 
$\ker(d_{\Hom}^A)$ is simply the graded endomorphisms of $C_A$ as complexes, and likewise for 
$C_B$. As $(C_A,\partial_A)$ is a dg-module over $(K,d_K)$ as well, and hence acyclic, 
we get that $\ker(d_{\Hom}^A)$ coincides with the graded endomorphism ring of $\ker(\partial_A)$  
over $\ker(d_K)$. Analogous statements hold for $C_B$. 
As $(C_A,\partial_A)$ is faithfully projective as a $K$-module, we get that 
$\ker(\partial_A)$ is faithfully projective as a graded $\ker(d_K)$-module.
Analogous statements hold for $\ker(\partial_B)$.

This shows that 
$$[\ker(d_A)]=[\ker(d_B)]$$
in $\grBr(\ker(d_K))$. 

\medskip

We need to show that the homomorphism has an inverse. 
In order to do so, let $[(A,d_A)]$ be an element in 
$\dgBr^I(K,d_K)$. 
Therefore the natural homomorphism
$$
\ker(d_A)\widehat\otimes_{\ker(d_K)}\ker(d_A)^{op}\lra \widehat\End_{\ker(d_K),graded}(\ker(d_A))
$$  
is an isomorphism of graded algebras. 
We also have the natural homomorphism of dg-algebras
$$
A\widehat\otimes_KA^{op}\lra \End_K^\bullet(A)
$$
which we need to show to be an isomorphism. 
However, as $(K,d_K)$ is acyclic, by Aldrich and Garcia-Rozas' theorem, we get that
also $(A,d_A)$ is acyclic, and again by the argument above,
$$K=\ker(d_K)\oplus\ker(d_K)T\mbox{ and }A=\ker(d_A)\oplus\ker(d_A)T$$
with the same $T$ in the (graded) centre of $A$. 
Hence,
\begin{eqnarray*}
A\widehat\otimes_KA^{op}&=&
(\ker(d_A)\widehat\otimes_{\ker(d_K)}K)\widehat\otimes_K
(\ker(d_A)\widehat\otimes_{\ker(d_K)}K)^{op}\\
&=&\left(\ker(d_A)\widehat\otimes_{\ker(d_K)}\ker(d_A)^{op}\right)\oplus
\left( \ker(d_A)\widehat\otimes_{\ker(d_K)}\ker(d_A)^{op}\right)T\\
&=&\widehat\End_{\ker(d_K),graded}(\ker(d_A))\oplus \widehat\End_{\ker(d_K),graded}(\ker(d_A))T\\
&=&\widehat\End_{K,graded}(A)
\end{eqnarray*}
using again that also $(\End_K^\bullet(A),d_\Hom)$ is a $(K,d_K)$-dg module and hence acyclic as well.
The map given by 
$$\dgBr^I(K,d_K)\ni [(A,d_A)]\mapsto [(A,d_A)]\in \dgBr^{{II}}(K,d_K)$$
is well-defined, by an argument similar as above. 
Indeed, if $[(A,d_A)]=[(B,d_B)]$ in $\dgBr^I((K,d_K))$, then there are faithfully projective 
graded $\ker(d_K)$-modules $P$ and $Q$ such that 
$$\ker(d_A)\widehat\otimes_{\ker(d_K)}\widehat\End_{\ker(d_K),graded}(P)\simeq 
\ker(d_B)\widehat\otimes_{\ker(d_K)}\widehat\End_{\ker(d_K),graded}(Q).$$
Put $(\widehat P,\partial_P):=P\widehat\otimes_{\ker(d_K)}(K,d_K)$ and 
$(\widehat Q,\partial_Q):=Q\widehat\otimes_{\ker(d_K)}(K,d_K)$.
Then 
\begin{eqnarray*}
(A,d_A)\widehat\otimes_K\End^\bullet_K(\widehat P,\partial_P)&=&
\ker(d_A)\widehat\otimes_{\ker(d_K)}(K,d_K)\widehat\otimes_K\End^\bullet_K(P\widehat\otimes_{\ker(d_K)}(K,d_K))\\
&=&\ker(d_A)\widehat\otimes_{\ker(d_K)}K\widehat\otimes_K(\widehat\End_{\ker(d_K),graded}(P)\widehat\otimes_{\ker(d_K)}K)\\
&=&\ker(d_A)\widehat\otimes_{\ker(d_K)}\widehat\End_{\ker(d_K),graded}(P)\widehat\otimes_{\ker(d_K)}K\\
&=&\ker(d_B)\widehat\otimes_{\ker(d_K)}\widehat\End_{\ker(d_K),graded}(Q)\widehat\otimes_{\ker(d_K)}K\\
&=&\ker(d_B)\widehat\otimes_{\ker(d_K)}K\widehat\otimes_K(\widehat\End_{\ker(d_K),graded}(Q)\widehat\otimes_{\ker(d_K)}K)\\
&=&\ker(d_B)\widehat\otimes_{\ker(d_K)}(K,d_K)\widehat\otimes_K\End^\bullet_K(Q\widehat\otimes_{\ker(d_K)}(K,d_K))\\
&=&(B,d_B)\widehat\otimes_K\End^\bullet_K(\widehat Q,\partial_Q)
\end{eqnarray*}
Also the fact that
this map is a group homomorphism follows by Lemma~\ref{Fcommuteswithtensor}. 
The fact that this map is left and right inverse to the map 
\begin{eqnarray*}
\dgBr^{II}(K,d_K)&\lra& \dgBr^{I}(K,d_K)\\
\ [(A,d)]&\mapsto&[(A,d)]
\end{eqnarray*}
from above is clear by definition. 
Hence the Theorem follows. \dickebox

\begin{Rem}
Lemma~\ref{dgBrIisgrBr} together with Proposition~\ref{dgBrIIvsusgradedBrK} show a partial analogous statement for dg-base rings with differential $0$. Here we treat the case of 
acyclic dg-base ring. In the case of differential $0$ we only get a splitting, with kernel
$L$ being the subgroup generated by the dg-algebras $A$, 
where $A$ is a graded algebra over $\ker(d_K)$ isomorphic to the graded endomorphisms 
of some faithfully projective 
$\ker(d_K)$-module, but with a non good grading. 

We may hence consider $L$ to be the group of bad gradings on 
matrix algebras (cf \cite{DascalescuIonNastasescuRioMontes}), 
modulo the subgroup of good gradings in an appropriate sense. 
\end{Rem}

We summarize our results up to now in the following scheme:

$$\xymatrix{
\dgBr^I(K,d_K)\ar[rrr]^\simeq_{\mbox{\scriptsize Lemma~\ref{dgBrIisgrBr}}}
&&&\grBr(\ker(d_K))\\ \\ \\ \\
\dgBr^{II}(K,d_K)\ar[uuuu]_\simeq^{\mbox{\scriptsize Theorem~\ref{theultimateforacyclic} for }(K,d_K)-dgmod \textup{ semisimple}}
\ar[rrr]_{\mbox{\scriptsize Theorem~\ref{dgBrIIvsusgradedBrK}}}
&&&\grBr(K)\ar@{_{(}->}@/_3pc/[lll]_{\mbox{\scriptsize split if }d_K=0, 
\mbox{\scriptsize Theorem~\ref{dgBrIIvsusgradedBrK}}}
}$$

\begin{Rem}
One may ask if there is a more natural definition of a dgBrauer group for 
a graded-commutative dg-base ring $(R,d_R)$. One could 
say that a dg-algebra $(A,d_A)$ over $(R,d_R)$ 
is dg-Azumaya of the third kind if 
\begin{itemize}
\item $(A,d_A)$ is a finitely generated projective object in the category of $(R,d_R)$-modules and
\item the natural homomorphism 
$(A,d_A)\widehat\otimes_R(A,d_A)^{op}\lra (\End_R^\bullet(A,d_A),d_{\Hom})$ is an isomorphism of dg-algebras. 
\end{itemize}
This definition parallels the classical definition of an Azumaya algebra naturally 
replacing each notion by the corresponding notion in the category of dg-modules, resp. dg-algebras.
However, we will get back $\dgBr^{II}(R,d_R)$.
Indeed, the proof of Lemma~\ref{groupstructureofdgbrauerofsecondkind} shows that 
the graded isomorphism $A\widehat\otimes_R A^{op}\lra \widehat\End_{R,graded}(A)$ gives rise to a 
dg-isomorphism 
$$(A\widehat\otimes_R A^{op},d_{A\widehat\otimes_R A^{op}})\lra (\End_{R}^\bullet(A,d_A),d_{\Hom}).$$ 

However, the hypothesis that $(A,d_A)$ should be a projective object in 
the category of $(R,d_R)$-modules is more restrictive. Indeed, this implies that 
$A$ is projective as a graded $R$-module, and in addition that $(A,d_A)$ is acyclic
(cf \cite[Proposition 3.3; proof of Proposition 3.4]{Tempest-Garcia-Rochas}). 
Since the unit element in such a dg-Brauer group would has to be acyclic again, 
we would have to impose that $(R,d_R)$ is acyclic, which then leads to the situation 
we already studied in Theorem~\ref{theultimateforacyclic}.

Note however that the graded centre of an acyclic dg-algebra 
$(A,\delta_A)$ is not necessarily acyclic. 
An example is given by the algebra 
$$\left(\begin{array}{cc}K&K\\ K&K\end{array}\right)=:A$$
with differential 
$$\delta_A(\left(\begin{array}{cc}a&b\\ c&d\end{array}\right))=
\left(\begin{array}{cc}c&d-a\\ 0&c\end{array}\right)$$
for a field $K$ and the obvious grading. 
This algebra is acyclic and the graded centre is $K\cdot 1_A$.
This latter algebra is not acyclic, and there is no acyclic dg-division 
algebra mapping onto $Z_{gr}(A)$.  
\end{Rem}

\section{Examples: DgBrauer groups of dg-fields}
\label{Examplesect}

\label{Examplesection}

Recall the different classes of dg-fields from Theorem~\ref{fieldclassscases}.
We shall now apply our results to each of these cases and use the same 
numbering as in the theorem.

Further, recall the following fact:
By \cite[IV.1.8 Theorem]{caenepeelvanoystaeyen} we get 
$$\ogrBr(L[T,T^{-1}])=\Br(L)\oplus H^2_{gr}(Gal(L^{sep}/L),|T|\Z),$$
where, denoting by $L^{sep}$ the separable closure of $L$,  
$$H^2_{gr}(Gal(L^{sep}/L),|T|\Z)=\ker(H^2(Gal(L^{sep}/L),|T|\Z)\lra H^2(Gal(L^{sep}/L),\Z)).$$
If the characteristic of $L$ does not divide $|T|$, then 
$$\ogrBr(L[T,T^{-1}])=\Br(L)\oplus H^2_{gr}(\Z/|T|\Z,(L^{sep})^\times)$$
A more precise behaviour depends on the question if $L$ is perfect or not
\cite[IV.1.10 Example]{caenepeelvanoystaeyen}. 

Recall  that for $K$ being a $\Z$-graded base ring 
concentrated in even degrees, we obtain a short exact sequence 
$$0\lra \ogrBr(K)\lra\grBr(K)\lra Q_2^{\Z}(K)\lra 0$$
of abelian groups, where $Q_2^{\Z}(K)$ is an abelian $2$-group of exponent at most $4$.
(cf \cite[Theorem 2.3]{TilborghsVanOystaeyen} and remarks at the beginning of Section 2. of loc.cit.)

\begin{enumerate}
\item{$K$ an ordinary field with zero differential}

Then Theorem~\ref{dgBrIIvsusgradedBrK} and Lemma~\ref{dgBrIisgrBr} show that 
$$\dgBr^{II}(K)=\grBr(K)=\dgBr^I(K).$$

\item{Laurent polynomials $K=L[T,T^{-1}]$ with zero differential}

Again, if $T$ is of even degree, 
Theorem~\ref{dgBrIIvsusgradedBrK} and Lemma~\ref{dgBrIisgrBr} show that 
$$\dgBr^{II}(L[T,T^{-1}])=\grBr(L[T,T^{-1}])=\dgBr^I(L[T,T^{-1}]).$$
If $T$ is not of even degree, the second equality still holds.

\item {Laurent polynomials $K=L[T,T^{-1}]$ with $d(T)=1$}.

Then $(K,d_K)$ is acyclic and its dg-module category is semisimple. Theorem~\ref{theultimateforacyclic} shows that 
$$\grBr(L[T^2,T^{-2}])=\dgBr^I(L[T,T^{-1}]).$$
If $L$ is of characteristic $2$, then commutativity and graded commutativity coincide, and 
$\dgBr^I(L[T,T^{-1}])=\dgBr^{II}(L[T,T^{-1}])$ here as well.

\item
\begin{enumerate}
\item{acyclic algebras $K[Y]/Y^2$ for some field $K$ and $d(Y)=1$}

As $K[Y]/Y^2$ is acyclic with $d(Y)=1$, so are all dg-algebras over $K[Y]/Y^2$.
Further the dg-module category over $K[Y]/Y^2$ is semisimple.
Theorem~\ref{theultimateforacyclic} then shows that
$$\dgBr^{II}(K[Y]/Y^2,d(Y)=1)=\dgBr^{I}(K[Y]/Y^2,d(Y)=1)=\grBr(K).$$

\item

{acyclic algebras $L[T,T^{-1}][Y]/Y^2$ for $d(Y)=1$}.
 
Again, as $L[T,T^{-1}][Y]/Y^2$ for $d(Y)=1$ is acyclic, so are all 
dg-algebras over the algebra $L[T,T^{-1}][Y]/Y^2$ for $d(Y)=1$.
Again the dg-module category of the base ring is semisimple.
Theorem~\ref{theultimateforacyclic} then shows that
\begin{eqnarray*}
\dgBr^{II}((L[T,T^{-1}][Y]/Y^2),d(Y)=1)&=&\dgBr^{I}((L[T,T^{-1}][Y]/Y^2),d(Y)=1)\\
&=&\grBr(L[T,T^{-1}]).
\end{eqnarray*}
\end{enumerate}

\item
As in the fifth case of the dg-field classification we needed to assume that the field is commutative, explicitly not graded-commutative, and in order to define dg-Brauer groups it is necessary to have a graded-commutative base-field, we cannot provide dg-Brauer groups in this case.
\end{enumerate}

\appendix

\section{Defining the graded Brauer group of a graded commutative ring}
\label{Appendixsection}

Let $R$ be a graded-commutative graded ring and $A$ a graded $R$-algebra. 
By definition, there is an injection from $R$ into the graded center of $A$. 
In this case, if we ignore the grading, $A$ is not necessarily an $R$-algebra. 
Hence, it is necessary to give a definition of a graded Azumaya algebra that 
does not rely on the ordinary notion of an Azumaya algebra. This is the purpose 
of this appendix.

\begin{Def} \cite{dgseparable}\label{dgseparabledef}
 A graded $R$-algebra $A$ is called a {\em graded separable algebra} 
 if it is projective as graded $A$-bimodule, or equivalently the multiplication map 
 $A\widehat\otimes_RA^{op}\lra A$ is a split epimorphism of graded $A\widehat\otimes_RA^{op}$-modules.
\end{Def}

\begin{Rem}
By \cite[Proposition 1.2.15]{Hazrat}, a graded $R$-module $P$ is a projective 
graded module if and only if $P$ is a graded $R$-module and $P$ is a projective 
$R$-module. Therefore, most properties are similar to those in the ungraded case.
\end{Rem}

For a graded $A$-bimodule $M$, denote by $M^A=\bigoplus_{n\in\Z} (M^{A})_n$, where 
$$(M^{A})_n:=\{m\in M_n|am=(-1)^{|a||m|}ma,\textup{ for any homogenous element }a\in A\}.$$ 
It is easy to see that $M^{A}$ isomorphic to the graded $R$-module $\widehat\Hom_{A^{e}}(A,M)$ and 
$A^A=Z_{gr}(A)$. A graded $R$-derivation $\partial:A\rightarrow M$ is an $R$-linear graded map 
satisfying: $\partial(ab)=\partial(a)b+(-1)^{|a||\partial|}a\partial(b)$. Denote by $J(A)$ 
the kernel of the morphism $A^e\rightarrow A$ as $A^e$-module, and by $\widehat{Der}(A,M)$ 
the set of graded $R$-derivations of $A$ with values in $M$.

\begin{Prop}\label{gr-seperable}
For the graded $R$-algebra $A$. The following conditions are equivalent:
\begin{enumerate}
\item
 $A$ is a graded separable $R$-algebra.
\item The functor $M\mapsto M^A$ is exact.
\item The $R$-derivation $\delta:A\lra A$ i.e. $\delta(a)=a\otimes 1-1\otimes a$ is an inner derivation.
\item All $R$-derivations are inner derivation.
\end{enumerate}
\end{Prop}
The proof similar to the ungraded case \cite[III, Theorem 1.4]{Knus},
using the isomorphism $M^A\cong \widehat\Hom_{A^{e}}(A,M)$, $A\delta(A)=J(A)$ 
and the isomorphism $\widehat\Hom_{A^{e}}(J(A),M)\cong \widehat{Der}(A,M)$.

\begin{Remark}\label{properties of graded seperable algebra}
   Using this proposition, one can easily prove that a graded separable algebra 
   still has the same properties as in \cite[III, 1]{Knus}.
\begin{enumerate}
\item Let $S_i$ for  $i\in\{1,2\}$ be graded commutative graded $R$-algebras and $A_i$ graded separable algebra over $S_i$ for  $i\in\{1,2\}$. 
Then $A_1\widehat\otimes_R A_2$ is graded separable algebra over $S_1 \widehat\otimes_R S_2$, and $Z_{gr}(A_1\widehat\otimes_R A_2)=Z_{gr}(A_1)\widehat\otimes Z_{gr}(A_2)$. In particular, 
$S_i\widehat\otimes_{R}A_i$ is graded separable $S_i$-algebra, 
and $Z_{gr}(S_i\widehat\otimes_{R}A_i)=S_i\widehat\otimes Z_{gr}(A_i), i=1,2$.
\item The graded center $Z_{gr}(A)$ of a graded separable $R$-algebra $A$ is a direct summand of $A$.
\item Let $A$ and $B$ be graded $R$-algebras. If $B$ is faithfully projective $R$-module and $A\widehat\otimes_RB$ be graded separable $R$-algebra, then $A$ is graded separable over $R$.\\
\item Let $\phi:A\rightarrow B$ be an epimorphism of graded $R$-algebras. If $A$ is a graded separable $R$-algebra, then so is $B$, and the graded center of $B$ is the image under $\phi$ of the graded center of $A$. 
\end{enumerate}
\end{Remark}

\begin{Lemma}\label{gr-center}
The graded center of a graded simple ring $R$ is a graded field. 
 \end{Lemma}

Proof.
Let $a\in Z_{gr}(R)$ be a nonzero homogenous element, then $Ra$ is a graded two-sided ideal of $R$, thus $Ra=R$. There exists a homogenous element $b$ such that $ba=1$. Similarly, there exists a homogenous element $c$ such that $ac=1$. Clearly, $b=bac=c$. For any homogenous element $x\in R$, $$bx=bxac=(-1)^{|x||a|}baxc=(-1)^{|x||b|}xb.$$ 
Hence,  $Z_{gr}(R)$ is a graded commutative graded division ring.
\dickebox

\begin{Prop}
Every graded separable algebra $A$ over a graded field $R$ is graded semisimple.
\end{Prop}

Proof. 
 Every graded $A$-module is projective graded $R$-module, since every graded module 
 over a graded division ring is free \cite[Proposition 4.6.1]{NastasescuVanOystaen}. 
 Moreover, for any graded $A$-modules $M, N$, we have an isomorphism 
 $\widehat\Hom_A(M,N)\cong \widehat\Hom_R(M,N)^A$. It follows that every graded $A$-module is 
 projective as graded $A$-module. In particular, every graded ideal of $A$ 
 is a direct summand of $A$. Therefore, it is easy to prove that $A$ is graded semisimple.   
\dickebox

\begin{Cor}
 Let $A$ be a graded separable algebra over a graded field. 
 If $A$ is graded central, then $A$ is graded simple algebra.
\end{Cor}

Proof.
In \cite[Section 2.9]{NastasescuVanOystaen}, 
$A=A_1\times \cdots \times A_n$, for some graded simple algebras 
$A_i$ for $i\in\{1,\cdots,n\}$. Since 
$$Z_{gr}(A_1\times \cdots \times A_n)=Z_{gr}(A_1)\times \cdots \times Z_{gr}(A_n)$$ 
and from Lemma~\ref{gr-center} it follows that $n=1$. Therefore $A$ is graded simple.   
\dickebox

\medskip

Recall the localization of graded rings \cite[Chapter 8]{NastasescuVanOystaen}. 
Suppose that $R$ is graded commutative, it is clear that every multiplicative 
closed subset $S$ of $R$ consisting of homogenous elements satisfies the left 
and right Ore conditions. Therefore, $S^{-1}R$ is a graded ring. 

\begin{Lemma}
 Let $\wp$ be a gr-prime of $R$ and $S$ a subset consisted of the homogenous elements of $R\setminus\wp$. The graded ring $R_\wp:=S^{-1}R$ is a gr-local ring.
\end{Lemma}

Proof.
For any $a,b \in S$, we have $a,b \notin \wp$. 
By the definition of a gr-prime ideal, it follows that $ab\notin \wp$. Therefore, 
$S$ is a multiplicative closed subset. Clearly, we can define a grading by
$$(\wp R_\wp)_n=\left\{\frac{p}{s}|s\in S, p\in \wp, |p|-|s|=n\right\},$$ 
and $\wp R_\wp$ is the unique maximal graded ideal. 
Indeed, let $I$ be a graded ideal containing $\wp R_\wp $. 
Suppose that there exists $\frac{r}{s}\in I\setminus \wp R_\wp$. Then $r\in R\setminus\wp$, 
and hence $\frac{s}{r}\in S^{-1}R$. Therefore, $1=\frac{s}{r}\frac{r}{s}\in I$. 
On the other hand, suppose that $I$ is a maximal graded ideal of $S^{-1}R$. 
If there exists $\frac{r}{s}\in I$ with $r\notin \wp$, then 
$1=\frac{s}{r}\frac{r}{s}\in I$, which is absurd. Therefore, 
$I\subseteq \wp R_\wp$, and since $I$ is maximal, it follows that $I=\wp R_\wp$.  
\dickebox

\begin{Lemma}\label{flatness is local property}
Let $A$ be a graded $R$-algebra and $M$ a graded $A$-module, 
then the following conditions are equivalent:
\begin{enumerate}
\item $M$ is graded flat $A$-module;
\item $M_m=A_m\otimes_AM=R_m\otimes_RA\otimes_AM$ is graded flat $A_m$-module 
for any maximal graded ideal $m$ of $R$.
\end{enumerate}
\end{Lemma}

Proof.
 (1) $\Rightarrow$ (2) is trivial.\\
 (2) $\Rightarrow$ (1) 
 Claim 1: $M=0$ if and only if $M_m=0$ for any maximal graded ideal $m$.\\
 $''\Rightarrow''$ is trivial. \\
 $''\Leftarrow''$  let $x\in M$ be a homogenous element, 
 $ann(x)=\{r\in R | rx=0\}$. It is easy to check that $ann(x)$ 
 is a graded ideal of $R$. For any maximal graded ideal $m$ we get by hypothesis
 $\frac{x}{1}=0$ in $M_m$, i.e. $\exists s\in R\setminus m$ such that $sx=0$. 
 Therefore, $s\in ann(x)\setminus m$, it follows that $ann(x)$ not contained in 
 any maximal graded ideal. Hence, $ann(x)=R$, which shows $M=0$.\\
 Claim 2: Suppose $L_1\xrightarrow{f} L_2 \xrightarrow{g} L_3$ is a sequence in the category of graded $R$-module, then it is exact if and only if $(L_1)_m\xrightarrow{f_m} (L_2)_m \xrightarrow{g_m} (L_3)_m$ is exact for any maximal graded ideal $m$.\\
 $'' \Rightarrow'' $ since $R_m$ is a graded flat $R$-module. \\
 $''\Leftarrow''$  $g_mf_m=(1\otimes g)(1\otimes f)=1\otimes gf=0$, which implies that 
 $\im (f)\subseteq \ker (g)$ since $m$ is was arbitrary.
$$(\ker (g)/ \im (f))_m \cong (\ker (g))_m/(\im (f))_m\cong \ker (g)_m/\im (f)_m=0,$$ for any maximal graded ideal $m$. Therefore, $\ker (g)/\im (f)=0,$ which implies $\im (f)= \ker (g)$.\\

We now finish the proof of (2) $\Rightarrow$ (1). 
Given a short exact sequence of graded right $A$-modules 
$$0\rightarrow L_1\rightarrow L_2\rightarrow L_3\rightarrow 0,$$ 
since $R_m$ is a graded flat $R$-module,  we have the exact 
sequence $$0\rightarrow (L_1)_m\rightarrow (L_2)_m\rightarrow (L_3)_m\rightarrow 0.$$ 
By hypothesis, $M_m$ is graded flat $A_m$-module, thus we have the exact sequence 
$$0\rightarrow (L_1)_m\otimes_{A_m}M_m\rightarrow (L_2)_m\otimes_{A_m}M_m\rightarrow (L_3)_m\otimes_{A_m}M_m\rightarrow 0$$
 It is clear that $(L_i)_m\otimes_{A_m}M_m\cong (L_i\otimes_AM)_m$ for every $i\in\{1,2,3\}$. 
 Therefore, for any maximal graded ideal $m$, we have the exact sequence
$$0\rightarrow  (L_1\otimes_AM)_m\rightarrow  (L_2\otimes_AM)_m\rightarrow  (L_3\otimes_AM)_m\rightarrow 0$$
 By claim 2, 
 $$0\rightarrow  L_1\otimes_AM\rightarrow L_2\otimes_AM\rightarrow  L_3\otimes_AM\rightarrow 0$$ is exact.
\dickebox

\begin{Prop}\label{propa10}
Let $A$ be a graded $R$-algebra, and it is finite generated as $R$-module. Then the following conditions are equivalent:
\begin{enumerate}
\item $A$ is graded separable $R$-algebra.
\item $A_\wp$ is graded separable $R_\wp$-algebra, for any gr-prime ideal $\wp$.
\item $A_m$ is graded separable $R_m$-algebra, for any maximal graded ideal $m$.
\end{enumerate}
\end{Prop}

Proof.
(1) $\Rightarrow$ (2)   $A_\wp=R_\wp\otimes_R A$, by remark \ref{properties of graded seperable algebra}, $A_\wp$ is graded separable algebra over $R_\wp$.\\
(2) $\Rightarrow$ (3) is trivial.\\
(3) $\Rightarrow$ (1) Since $A_m$ is graded separable over $R_m$, that is 
$A_m$ is graded projective $(A_m)^e$-module. Thus, $A_m$ is graded flat 
$(A_m)^e$-module. Therefore, by the lemma \ref{flatness is local property}, 
$A$ is graded flat. Since $A$ is finite generated $R$-module and 
$J(A)=A\delta(A)$, then it is clear that $J(A)$ is also a finite generated $R$-module. 
It follows that $A$ has a finite presentation as $A^e$-module.\\
Now use the following two facts:
\begin{enumerate}
\item
 A module of finite presentation is flat if and only if it is a finitely generated projective;
\item A graded module $M$ is a graded flat module over a graded ring $R$ if and only if $M$ is 
a graded and flat $R$-module \cite[Proposition 3.4.5]{Hazrat}.
\end{enumerate}
Therefore, $A$ is a projective $A^e$-module. 
By \cite[Proposition 1.2.15]{Hazrat}, $A$ is graded projective $A^e$-module, 
whence graded-separable over $R$.
\dickebox

\begin{Prop}\label{maximal graded ideal}
Let $A$ be a graded $R$-algebra which is finitely generated as $R$-module.
Then $A$ is graded separable 
 if and only if $A/mA$ is graded separable over $R/m$, for any maximal graded ideal $m$.   
\end{Prop}

Proof.
$''\Rightarrow''$  By the remark \ref{properties of graded seperable algebra}, 
$A/mA=A\otimes_RR/m$ is graded separable over $R/m$, for any maximal graded ideal $m$.\\
$''\Leftarrow''$ It is easy to prove that $R_m/mR_m\cong R/m$. Then we have 
$$A_m/mA_m=(A\otimes_RR_m)/m(A\otimes_RR_m)\cong A\otimes_R(R_m/mR_m)\cong A\otimes_RR/m\cong A/mA.$$ 
Therefore, we can prove the proposition for a graded local ring $R$ first. 
Here the proof is similar as the ungraded case \cite[Proposition 2.6]{Knus}. 
Using proposition \ref{gr-seperable} and the graded Nakayama's lemma then gives the first step. Then, use Proposition~\ref{propa10}.
\dickebox

\begin{Def}\label{gradedAzumayadefAppendix}
Let $R$ be a graded commutative graded ring.  
A graded $R$-algebra $A$ is called a {\em graded Azumaya algebra} 
if it is both graded central and graded separable.
\end{Def}

\begin{Prop}\label{graded Azumaya}
Let $R$ be a graded-commutative graded ring and $A$ a graded $R$-algebra. The following conditions are equivalent:
\begin{enumerate}
\item\label{(1)} $A$ is graded Azumaya algebra over $R$.
\item\label{(2)} $A$ is faithfully projective graded $R$-module, $a\otimes b\mapsto(x\mapsto (-1)^{|x||b|}axb)$ induces an isomorphism of graded algebras $A\widehat\otimes_RA^{op}\lra \widehat\End_{R,graded}(A)$.
\item\label{(3)} The functor $N\mapsto A\widehat\otimes_R N$ and $M\mapsto M^A$ establish a 
graded equivalence of categories of categories of graded $R$-modules and graded $A^e$-modules.
\item\label{(4)} $A$ is finitely generated graded $R$-module, $A/mA$ is a 
graded Azumaya algebra over $R/m$, 
for any maximal graded ideal $m$ of $R$.
\item\label{(5)} There exists a graded $R$-algebra $B$ and a faithfully 
projective graded $R$-module $P$ such that $A\widehat\otimes_RB\lra \widehat\End_{R,graded}(P)$ 
as graded $R$-algebras.
\end{enumerate}
\end{Prop}

Proof.
(\ref{(1)})$\Rightarrow$ (\ref{(2)})
By \cite[Proposition 1.2.15]{Hazrat} and \cite[Proposition 2.2.5]{Hazrat},
there exists a Morita theory analogous to that in \cite[I.7]{Knus}. 
Using these propositions and lemmas, together with the graded version 
of the dual basis lemma \cite[Theorem 1.2.17]{Hazrat}, one can follow 
the proof in the ungraded case given in \cite{Knus}. 
 
(\ref{(2)})$\Rightarrow$ (\ref{(3)}), (\ref{(3)})$\Rightarrow$ (\ref{(1)}), 
(\ref{(2)})$\Rightarrow$(\ref{(5)}) and (\ref{(5)})$\Rightarrow$(\ref{(1)}) 
are proved similarly to the ungraded case \cite{Knus}, 
using general graded Morita theory \cite[2.3]{Hazrat}. 
 
(\ref{(1)})$\Rightarrow$ (\ref{(4)}) follows the proof of (\ref{(1)})$\Rightarrow$ (\ref{(2)}).
 
(\ref{(4)})$\Rightarrow$ (\ref{(1)}) By the proposition~\ref{maximal graded ideal}, 
$A$ is a 
graded separable $R$-algebra. By the proposition~\ref{properties of graded seperable algebra}, 
$Z_{gr}(A/mA)=Z_{gr}(A)/Z_{gr}(A)\cap mA=Z_{gr}(A)/mZ_{gr}(A)$. Since $A/mA$ is graded central 
over $R/m$, thus $Z_{gr}(A/mA)=R/m$, it follows that $R/m\cong Z_{gr}(A)/mZ_{gr}(A)$, for any 
maximal graded ideal $m$. Then, similar as the proof of lemma \ref{flatness is local property}, 
one can prove that $R\cong Z_{gr}(A)$.
\dickebox

\medskip

Consider an equivalence relation $''\sim\ ''$ on graded $R$-algebras: $A\sim B$ 
if and only if there exist faithfully projective graded $R$-modules $P$ and $Q$ 
such that the graded ring $A\widehat\otimes_R\widehat\End(P)$ is isomorphic to the graded ring 
$B\widehat\otimes_R\widehat\End(Q)$. One 
can prove that this is indeed an equivalence relation.\\
Let $\grBr(R)$ denote the set of equivalence class of all graded Azumaya 
$R$-algebras with respect to the equivalence relation $''\sim\ ''$. 
By the remark \ref{properties of graded seperable algebra}, the operation of tensor 
product over $R$ is compatible with the equivalence relation so that there is induced 
an associative and commutative multiplication in $\grBr(R)$. The equivalence class which 
contains $R$ itself is clearly an identity for this multiplication.  If $A$ is a graded 
Azumaya $R$-algebra, then clearly $A^{op}$ is also a graded Azumaya $R$-algebra. 
And by the proposition \ref{graded Azumaya}, we have that $A$ is a faithfully 
projective graded $R$-module and that $A\widehat\otimes_RA^{op}=\widehat\End_R(A, A)$. 
Therefore it follows that $A\widehat\otimes_RA^{op}\sim R$ so that the equivalence 
class of $A^{op}$ is an inverse to that of $A$ in $\grBr(R)$. Thus, we have 
proved that $\grBr(R)$ is a group.

\medskip

Recall that for a graded $R$-algebra $B$ and a graded $B$-bimodule $M$ we denote by
$M^B$ the module with degree $n$ homogeneous component 
$\{m\in M_n\;|\;bm=(-1)^{n|b|}mb\;\forall b\textup{ homogeneous}\}$.

\begin{Prop} (the graded-double commutant theorem) \label{graded-double-commutant}
Let $R$ be a graded-commutative algebra and let $C$ be a graded $R$-algebra with 
$Z_{gr}(C)=R$. Let $A\subseteq C$ be a graded subalgebra
with $Z_{gr}(A)=C$, and such that $A$ is separable over $R$. Then 
\begin{enumerate}
\item the map
$A\widehat\otimes_RC^A\lra C$ given by $a\otimes c\mapsto ac$ is an isomorphism.
\item
$Z_{gr}(C^A)=R$.
\item
If $C$ is graded separable over $R$, then $C^A$ is graded-separable over $R$.
and  $C^{(C^A)}=A$.
\end{enumerate}
\end{Prop}

Proof. 
Since $A$ is graded-separable over $R$, 
by Proposition~\ref{graded Azumaya}.\ref{(1)}
we also get that $A$ is faithfully projective over $R$.
Then, the evaluation map $A\widehat\otimes_R\widehat\Hom_{\widehat\End_R(A)}(A,C)\lra C$
is an isomorphism. Since $A$ is graded-separable over $R=Z_{gr}(A)$, we have  
$\widehat\End_R(A)\simeq A\widehat\otimes_RA^{op}$, and since 
$\Hom_{A\widehat\otimes_RA^{op}}(A,C)=C^A$, we observe that then the evaluation map is 
an algebra homomorphism, which proves the first statement. 

Let $x\in Z_{gr}(C^A)$ be homogeneous. Then, since $A\widehat\otimes_RC^A\stackrel\mu\lra C$ 
is an isomorphism, we have for homogeneous $a\in A$ and $c\in C^A$ that
$$\mu(a\otimes c)x=acx=(-1)^{|c||x|}axc=(-1)^{(|c|+|a|)|x|}xac=(-1)^{(|c|+|a|)|x|}x\mu(a\otimes c),$$ 
hence $x$ graded-commutes with $C$. 
Hence, $Z_{gr}(C^A)\subseteq R$. The other inclusion is clear.

Since $A$ is separable over $R$, we get that $A$ is finitely generated 
graded faithfully projective over $R$, and $R$ is a direct summand 
of $A$ as an $R$-module, $A\otimes_RC^A\simeq C$ and $C$ separable over $R$ 
implies that $C^A$ is separable over $R$.

Let $B=C^{(C^A)}$. Then $A\subseteq B\subseteq C$. Since we already know that 
$C^A$ is graded-separable over $R$, we see that the multiplication map is an isomorphism
$$C^A\widehat\otimes_RB\simeq C.$$
But also the multiplication map $$A\widehat\otimes_RC^A\lra C$$ is an isomorphism.
Since $C^A$ is graded-separable over $R$, also $C^A$ is graded-projective over $R$, and hence 
$$(B/A)\widehat\otimes_RC^A=(B\widehat\otimes_RC^A)/(A\widehat\otimes_RC^A)=0.$$
This shows that $B=A$, using again that $C^A$ is graded-faithfully projective over $R$. 
\dickebox

\begin{Cor}\label{centralizersoftensorprodofAzumayas}
Let $R$ be a graded-commutative ring, and let $A$ and $B$ be graded-Azumaya algebras.
Then $(A\widehat\otimes_RB)^B=A$ and $(A\widehat\otimes_RB)^A=B$.
\end{Cor}

Proof. Use Proposition~\ref{graded-double-commutant}  for $C=A\widehat\otimes_RB$. 
Note that 
$A\widehat\otimes_RB$ contains a copy $A\widehat\otimes_R1$
of $A$ and a copy $1\widehat\otimes_RB$ of $B$.
We have that multiplication on factors yield
$A\widehat\otimes_RC^A\simeq C\simeq A\widehat\otimes_RB$, and 
since $B\subseteq C^A$, we get that $A\widehat\otimes_R(C^A/B)=0$.
Since $A$ is faithfully projective over $R$, this shows $B=C^A$.
\dickebox

\end{document}